\documentclass[numbers,webpdf,imaiai]{ima-authoring-template}%

\newcommand{\f}{\mathbf}

\newcommand{\lr}[1]{\left\langle{#1}\right\rangle}
\newcommand{\krylov}[3]{\mathcal{K}_{#1}({#2}, {#3})}

\newcommand{\Lg}{L_{\f{g}}}
\newcommand{\inH}{\f{\bar{H}}}
\newcommand{\eH}{\f{H}}
\newcommand{\sigmoid}{\text{sigmoid}}

\renewcommand\th{\textsuperscript{th}}
\newcommand*\dotprod[1]{\left\langle #1\right\rangle}
\newcommand*\vnorm[1]{\left\| #1\right\|}

\DeclareMathOperator*{\argmin}{arg\,min}

\usepackage[capitalise]{cleveref}
\crefname{equation}{}{}
\crefname{figure}{Figure}{Figures}
\creflabelformat{equation}{\textup{(#2#1#3)}}
\crefname{assumption}{Assumption}{Assumptions}
\crefname{condition}{Condition}{Conditions}
\crefname{corollary}{Corollary}{Corollaries}
\crefname{lemma}{Lemma}{Lemmas}

\usepackage{dsfont, comment, amsmath, subcaption}

\graphicspath{{Fig/}}

\theoremstyle{thmstyletwo}%
\newtheorem{theorem}{Theorem}
\newtheorem*{theorem*}{Theorem}
\newtheorem{assumption}{Assumption}
\newtheorem{condition}{Condition}
\newtheorem{corollary}{Corollary}
\newtheorem{lemma}{Lemma}
\newtheorem{example}{Example}%
\newtheorem{remark}{Remark}%
\newtheorem{property}{Properties}
\newtheorem{definition}{Definition}

\numberwithin{equation}{section}

\begin{document}

\DOI{DOI HERE}
\copyrightyear{2021}
\vol{00}
\pubyear{2021}
\access{Advance Access Publication Date: Day Month Year}
\appnotes{Paper}
\copyrightstatement{Published by Oxford University Press on behalf of the Institute of Mathematics and its Applications. All rights reserved.}
\firstpage{1}


\title[Nonconvex Newton-MR Under Inexact Hessian]{Complexity Guarantees for Nonconvex Newton-MR Under Inexact Hessian Information}

\author{Alexander Lim
\address{\orgdiv{School of Mathematics and Physics}, \orgname{University of Queensland}, \orgaddress{\street{St Lucia}, \postcode{4072}, \state{Queensland}, \country{Australia}}}}
\author{Fred Roosta 
\address{\orgdiv{School of Mathematics and Physics}, \orgname{University of Queensland}, \orgaddress{\street{St Lucia}, \postcode{4072}, \state{Queensland}, \country{Australia}}}}

\authormark{Lim. A and Roosta. F}

\corresp[**]{Corresponding author: \href{fred.roosta@uq.edu.au}{fred.roosta@uq.edu.au}}

\received{Date}{0}{Year}
\revised{Date}{0}{Year}
\accepted{Date}{0}{Year}


\abstract{We consider an extension of the Newton-MR algorithm for nonconvex unconstrained optimization to the settings where Hessian information is approximated. Under a particular noise model on the Hessian matrix, we investigate the iteration and operation complexities of this variant to achieve appropriate sub-optimality criteria in several nonconvex settings. We do this by first considering functions that satisfy the (generalized) Polyak-\L ojasiewicz condition, a special sub-class of nonconvex functions. We show that, under certain conditions, our algorithm achieves global linear convergence rate. We then consider more general nonconvex settings where the rate to obtain first order sub-optimality is shown to be sub-linear. In all these settings, we show that our algorithm converges regardless of the degree of approximation of the Hessian as well as the accuracy of the solution to the sub-problem. Finally, we compare the performance of our algorithm with several alternatives on a few machine learning problems.}

\maketitle

\section{Introduction}
Consider the following unconstrained optimization problem:
\begin{align}\label{eq:min_f}
    \min_{\f{x}\in\mathbb{R}^d}f(\f{x}),
\end{align}
where $f : \mathbb{R}^d \to \mathbb{R}$ is twice continuously differentiable and nonconvex.
Modern machine learning, in particular deep learning, has motivated the development of a plethora of optimization methods for solving \cref{eq:min_f}. The majority of these methods belong to the class of first order algorithms, which use only gradient information \cite{beck2017first,tian2023recent,ruder2016overview,alfarra2020adaptive,goodfellow2016deep,lin2020accelerated,lan2020first,zeiler2012adadelta,duchi2011adaptive,kingma2014adam,loshchilov2017decoupled}.
While first order methods are typically memory efficient, have low-cost per iteration, and are simple to implement, they suffer from serious drawbacks.  These methods are notoriously difficult to fine-tune and tend to converge slowly, especially when dealing with ill-conditioned problems. 
By incorporating the curvature information in the form of (approximate) Hessian matrix, second order algorithms  \cite{cartis2022evaluation,nocedal2006numerical} tend to alleviate these drawbacks. Beyond their theoretical appeal such as affine invariance and fast local convergence, these methods are also empirically shown to be less sensitive to hyper-parameters tuning and problem ill-conditioning \cite{berahas2020investigation, xu2020second}. 
 
Arguably, the primary computational bottleneck in second order methods for solving large-scale problems is due to operations involving the Hessian matrix. In such large-scale settings where storing the Hessian matrix explicitly is prohibitive, operations such as matrix-vector products would impose a cost equivalent to multiple function evaluations; see  \cite{liu2021convergence,roosta2022newton}.
This is often considered as the major disadvantage of these methods. 
Consequently, second order algorithms that utilize appropriate approximations of the Hessian matrix have been introduced, e.g., methods based on constructing probabilistic models \cite{rinaldi2024stochastic,jin2024sample,bellavia2022stochastic,blanchet2019convergence,byrd2011use,gratton2018complexity,tripuraneni2018stochastic,li2023randomized}, methods that employ randomized sub-sampling in finite sum settings \cite{roosta2019sub,xu2016sub,xu2020newton,yao2021inexact,yao2023inexact,byrd2012sample,bergou2022subsampling,bellavia2020subsampled,bellavia2021adaptive,bellavia2020inexact,bellavia2022linesearch,bollapragada2016exact,berahas2020investigation,na2023hessian,chen2022san}, similar techniques based on randomized sketching of the Hessian matrix  \cite{berahas2020investigation,na2023hessian,pilanci2017newton,feng2021non,na2024statistic}, and the related randomized subspace methods that are at times categorized under the sketch-and-project framework  \cite{fuji2022randomized,cartis2022randomised,gower2019rsn,na2023hessian,yuan2022sketched,derezinski2024sharp}.
A bird’s eye view of these methods reveals three main categories of globalization strategies based on trust region \cite{conn2000trust,xu2020newton}, cubic regularization \cite{cartis2011adaptive,cartis2011adaptive2,nesterov2006cubic}, and line-search \cite{nocedal2006numerical}. 
Almost all second order methods with line-search algorithmic framework employ the conjugate gradient (CG) algorithm as their respective sub-problem solver, e.g., variants of Newton-CG \cite{roosta2019sub,bollapragada2016exact,yao2023inexact,bellavia2022linesearch,bellavia2020subsampled}. An exception lies in  \cite{liu2021convergence} where the Newton-MR algorithm, initially proposed in \cite{roosta2022newton}, is extended to incorporate Hessian approximation. Unlike Newton-CG, the framework of Newton-MR relies on the minimum residual (MINRES) method \cite{liu2022minres,paige1975solution} as the sub-problem solver. However, the Newton-MR variants discussed in \cite{liu2021convergence} are limited in their scope to invex problems \cite{mishra2008invexity}. 

In this paper, we extend the analysis in \cite{liu2021convergence} to accommodate inexact Hessian information in a variety of more general nonconvex settings. 
In addition to iteration complexity, we provide operation complexity of our method. 
This is so since in large-scale problems, where solving the sub-problem constitutes the main computational costs, analyzing operation complexity -- which incorporates the cost of solving these sub-problems -- better captures the essence of the algorithmic cost than simple iteration complexity \cite{cartis2022evaluation}. 
Although extensively studied in convex settings, efforts to establish operation complexity results for nonconvex settings have emerged only recently. For example, to our knowledge, \cite{carmon2018analysis} was among the first to analyze the complexity of solving sub-problems in trust region and cubic regularization methods in nonconvex settings, while \cite{royer2018complexity,royer2020newton} focused on line-search-based methods. Since then, numerous follow-up works have aimed to quantify the overall operational complexity across a wide range of nonconvex optimization algorithms for solving \cref{eq:min_f} by assessing the computational cost of solving their respective sub-problems, e.g., \cite{yao2023inexact,liu2022newton,arjevani2020second,li2023randomized,curtis2021trust}.

\textit{The overarching objective of our work here is to derive both iteration and operation complexities of our proposed nonconvex Newton-MR variant in the presence of inexact Hessian information with minimal assumptions and limited algorithmic modifications compared with similar Newton-type methods in the literature}. In particular, (i) we only consider minimal smoothness assumptions on the gradient and do not extend smoothness assumptions to the Hessian itself; and (ii) we do not employ the Hessian regularization/damping techniques used in \cite{yao2023inexact, royer2020newton, li2023randomized, curtis2021trust}, as distorting the curvature information further could potentially lead to inferior performance (see experiments in \cref{sec:hessiam_damp}). 

The rest of this paper is organized as follows. We end this section by presenting the essential notation and definitions used in this paper. In \cref{sec:newton_mr_review}, we lay out the details of our algorithm. This is then followed by the  theoretical analyses in \cref{sec:theoretical_analyses}. Specifically, we first establish the iteration complexities of our algorithm is \cref{sec:iteration}. Subsequently, we derive the underlying operation complexities in \cref{sec:operation}. In \cref{sec:numerical_exp}, we empirically evaluate the performance of our algorithm on several machine learning problems. Conclusions and further thoughts are gathered in \cref{sec:conclusion}. Some proof details as well as further numerical results are deferred to \cref{sec:appendix}.

\subsection*{Notation}
Throughout the paper, vectors and matrices are, respectively, denoted by bold lower and upper case letters, e.g. $\f{a}$ and $\f{A}$. 
Their respective norms, e.g., $\|\f{a}\|$ and $\|\f{A}\|$, are Euclidean and matrix spectral norms. Regular letters represent scalar, e.g. $d$, $L$, $\alpha$, etc. 
We denote the exact and the inexact Hessian as $\eH$ and $\inH$, respectively. We use subscripts to denote the iteration counter of our algorithm, e.g., $\f{x}_k$. 
The objective function and its gradient at iteration $k$ are denoted by $f_k \triangleq f(\f{x}_k)$ and $\f{g}_k \triangleq \f{g}(\f{x}_k)$, respectively. We use superscripts to denote the iterations of MINRES. 
Specifically, $\f{s}_k^{(t)}$ refers to the $t\th$ iterate of MINRES using $\inH_k$ and $\f{g}_k$ with the corresponding residual $\f{r}_k^{(t)} = -\inH_k\f{s}_k^{(t)} - \f{g}_k$. 
The Krylov subspace of degree $t$, generated by $\inH_k$ and $\f{g}_k$, is denoted by $\mathcal{K}_t(\inH_k, \f{g}_k)$. Natural logarithm is denoted by $\log(\cdot)$. 
The range and the nullspace of a matrix, say $\inH$, are denoted by $\text{Range}(\inH)$ and $\text{Null}(\inH)$, respectively.
We also denote $f^{\star}\triangleq\min f(\f{x})$, which is assumed finite.

\subsection*{Definitions}
We analyze the complexity of our proposed algorithm in nonconvex settings. An interesting subclass of nonconvex functions that often times allows for favourable convergence rates is those satisfying (generalized) Polyak-\L ojasiewicz (PL) condition.   
\begin{definition}[\texorpdfstring{$\theta$}{theta}-Polyak-\L ojasiewicz  Condition]\label{def:PL}
For any $\theta \in [1, 2]$, we say a function satisfies the $ \theta $-PL condition, if there exists a constant $ \mu > 0 $ such that 
\begin{align}
\label{cond:PL}
	\|\f{g}(\f{x})\|^{2} \geq 2\mu(f(\f{x}) - f^{\star})^{{2}/{\theta}}, \quad \forall \f{x} \in \mathbb{R}^{d}.
\end{align}
\end{definition}
The class of $\theta$-PL functions has been extensively considered in the literature for the convergence analysis of various optimization algorithms, e.g., \cite{karimi2016linear,roosta2022newton, mei2020global,agarwal2021theory,yuan2022general,allen2019convergence,zeng2018global,bassily2018exponential,fatkhullin2022sharp}. It has been shown that \cref{cond:PL} is a special case of global Kurdyka-\L ojasiewicz inequality \cite{fatkhullin2022sharp}. Note that as $\theta$ gets closer to $1$, the function is allowed to be more flat near the set of optimal points, e.g., consider $f(x) = x^{2\alpha}$ for $\alpha \geq 1$, which satisfies \cref{cond:PL} with $\theta = 2\alpha/(2\alpha-1)$ (and  $\mu = 0.5 (2\alpha)^{(2\alpha-1)/\alpha}$). 

In optimization, one's objective is often to find solutions that satisfy a certain level of approximate optimality. When dealing with $\theta$-PL functions, it is appropriate to consider approximate global optimality.
\begin{definition}[$\varepsilon_f$-Global Optimality]\label{def:global_sub-optimality}
Given $0 < \varepsilon_{f} < 1$, a point $\f{x} \in \mathbb{R}^d$ is an $\varepsilon_{f}$- global optimal solution to the problem \cref{eq:min_f} if
\begin{equation}\label{eq:global_optimality_condition}
    f(\f{x}) - f^* \leq \varepsilon_{f}
\end{equation}
\end{definition}  
For more general nonconvex functions, a more suitable measure of sub-optimality is that of the approximate first order criticality. 
\begin{definition}[$\varepsilon_\f{g}$-First Order Optimality]\label{def:1st-order_sub-optimality}
Given $0 < \varepsilon_{\f{g}} < 1$, a point $\f{x} \in \mathbb{R}^d$ is an $\varepsilon_{\f{g}}$- first order optimal solution to the problem \cref{eq:min_f} if
\begin{equation}\label{eq:1st_order_optimality_condition}
    \|\f{g}(\f{x})\| \leq \varepsilon_{\f{g}}
\end{equation}
\end{definition}  

The geometry of the optimization landscape is essentially encoded in the spectrum of the Hessian matrix. In that light, one can obtain descent by leveraging directions that, to a great degree, align with the eigenspace associated with the small or negative eigenvalues of the Hessian. 
These directions are referred to as $\sigma$-Limited Curvature (LC) directions. 
\begin{definition}[$\sigma$-Limited Curvature Direction]\label{def:LC}
For any nonzero $\f{x}\in\mathbb{R}^d$ and $\sigma > 0$, we day $\f{x}$ is a $\sigma$-limited curvature ($\sigma$-LC) direction for a  matrix $\f{A}$, if $\langle\f{x}, \f{A} \f{x}\rangle \leq d \sigma \|\f{x}\|^2.$
\end{definition}

We perform our convergence analysis under a simple noise model for the Hessian matrix. 
\begin{definition}[Inexact Hessian] \label{def:noise_model}
For some $\varepsilon > 0$, the matrix $\inH(\f{x})$ is an estimate of the underlying Hessian, $\eH(\f{x})$, in that 
\begin{align}\label{eq:noise_model}
    \vnorm{\eH(\f{x}) - \inH(\f{x})} \leq \varepsilon.
\end{align}
\end{definition}  
This type of inexactness model includes popular frameworks such as the sub-sampled Newton \cite{roosta2019sub,xu2020newton} and Newton Sketch \cite{pilanci2017newton,feng2021non} and is considered frequently in the literature, e.g., \cite{yao2023inexact,liu2021convergence,li2023randomized,bollapragada2016exact}. For example, consider large-scale finite sum minimization problems \cite{shalev2014understanding}, where $f(\f{x}) = \sum_{i=1}^n f_i(\f{x})/n$ and $n \gg 1$. Note that $\eH(\f{x}) = \sum_{i = 1}^n \f{H}_i(\f{x})/n$. Under a simple Lipschitz gradient assumption, it is shown in \cite[Lemma 16]{xu2020newton} that, for any $\varepsilon > 0$ and $0 < \delta < 1$, if $\inH(\f{x}) = \sum_{i \in \mathcal{S}} \f{H}_i(\f{x})/|\mathcal{S}|$ is constructed using $|\mathcal{S}| \in  \mathcal{O}(\varepsilon^{-2} \log (2d/\delta))$ samples drawn uniformly at random from $\{1,2,\ldots,n\}$, we obtain \cref{eq:noise_model} with probability at least $1-\delta$.

\section{Newton-MR with Inexact Hessian}\label{sec:newton_mr_review}
We now present our variant of Newton-MR that incorporates inexact Hessian information. The core of our algorithm lies on MINRES, depicted in \cref{alg:MINRES}, as the sub-problem solver. \cref{alg:MINRES} is similar to Algorithm 2.1 in \cite{liu2022minres}. The main differences are that \cref{alg:MINRES} is equipped with early termination criteria, outlined in \hyperref[alg:MINRES:inexactness]{Step 7} and \hyperref[alg:MINRES:inexactness]{Step 11}, which respectively correspond to \emph{sub-problem inexactness} (\cref{cond:inexactness}) and \emph{detection of limited curvature} (\cref{cond:LC}).

\begin{algorithm}[htbp]
	\caption{MINRES$(\inH, \f{g}, \eta, \sigma)$}
	\label{alg:MINRES}
	\begin{algorithmic}[1]
		\Require Inexact Hessian $\inH$, Gradient $\f{g}$, Inexactness tolerance $\eta > 0$, limited curvature tolerance $\sigma > 0$
		\State $\phi_0 = \tilde{\beta}_1 = \|\f{g}\|$, $\f{r}_0 = -\f{g}$, $\f{v}_1 = \f{r}_0/\phi_0$, $\f{v}_0 = \f{s}_0 = \f{w}_0 = \f{w}_{-1} = \f{0}$,
        \State $c_0 = -1$, $\delta_1^{(1)} = s_0 = \tau_0 = 0$, \texttt{Type = SOL},
        \State $t = 1$
		\While{True}
		\State $\f{q}_t = \inH\f{v}_t$, $\tilde{\alpha}_t = \f{v}_t^\top\f{q}_t$, $\f{q}_t = \f{q}_t - \tilde{\beta}_t\f{v}_{t-1} - \tilde{\alpha}_t\f{v}_t$, $\tilde{\beta}_{t+1} = \|\f{q}_t\|$
        \State $\begin{bmatrix}
            \delta_t^{(2)} & \epsilon_{t+1}\\
            \gamma_t^{(1)} & \delta_{t+1}^{(1)}
        \end{bmatrix} = \begin{bmatrix}
            c_{t-1} & s_{t-1}\\
            s_{t-1} & -c_{t-1}
        \end{bmatrix}\begin{bmatrix}
            \delta_t^{(1)} & 0 \\
            \tilde{\alpha}_t & \tilde{\beta}_{t+1}
        \end{bmatrix}$
        \If{$\phi_{t-1}\sqrt{\left(\gamma_t^{(1)}\right)^2 + \left(\delta_{t+1}^{(1)}\right)^2} \leq \eta\sqrt{\phi_0^2 - \phi_{t-1}^2}$} \label{alg:MINRES:inexactness} \Comment{\cref{cond:inexactness}: $\vnorm{\inH\f{r}_{t-1}} \leq \eta \vnorm{\inH\f{s}_{t-1}}$}
        \State \texttt{Type = SOL} \Comment{\texttt{SOL} direction}
        \State \Return $\f{s}_{t-1}$, \texttt{Type}
		\EndIf
        \If{$-c_{t-1}\gamma_t^{(1)} \leq \sigma d$} \label{alg:MINRES:LC} \Comment{\cref{cond:LC}: $\lr{\f{r}_{t-1},\inH\f{r}_{t-1}} \leq \sigma d \vnorm{\f{r}_{t-1}}^2$}
        \State \texttt{Type = LC} \Comment{\texttt{$\sigma$-LC} direction}
        \State \Return $\f{r}_{t-1}$, \texttt{Type}
        \EndIf
        \State $\gamma_t^{(2)} = \sqrt{\left(\gamma_t^{(1)}\right)^2 + \tilde{\beta}_{t+1}^2}$
        \If{$\gamma_t^{(2)} \not= 0$}
        \State $c_t = \gamma_t^{(1)}/\gamma_t^{(2)}$, $s_t = \tilde{\beta}_{t+1}/\gamma_t^{(2)}$, $\tau_t = c_t \phi_{t-1}$, $\phi_t = s_t \phi_{t-1}$,
        \State $\f{w}_t = \left(\f{v}_t - \delta_t^{(2)}\f{w}_{t-1} - \epsilon_t\f{w}_{t-2}\right) / \gamma_t^{(2)}$, $\f{s}_t = \f{s}_{t-1} + \tau_t \f{w}_t$
        \If{$\tilde{\beta}_{t+1}\not= 0$}
        \State $\f{v}_{t+1} = \f{q}_t / \tilde{\beta}_{t+1}, \f{r}_t = s_t^2\f{r}_{t-1} - \phi_t c_t \f{v}_{t+1}$
        \Else
        \State \Return $\f{s}_t$, \texttt{Type}
        \EndIf
        \Else
        \State $c_t = 0$, $s_t = 1$, $\tau_t = 0$, $\phi_t = \phi_{t-1}$, $\f{r}_t = \f{r}_{t-1}$, $\f{s}_t = \f{s}_{t-1}$
        \EndIf
        \State t = t + 1
		\EndWhile
	\end{algorithmic}
\end{algorithm}

\begin{condition}[Sub-problem Inexactness Condition]\label{cond:inexactness}
    If $\|\inH_k\f{r}_k^{(t-1)}\| < \eta\|\inH_k\f{s}_k^{(t-1)}\|$ at iteration $ t $ of \cref{alg:MINRES}, 
    then $ \f{s}_k^{(t-1)} $ is declared as a solution direction satisfying the sub-problem inexactness condition where $\eta > 0$ is any given inexactness parameter. This condition is effortlessly verified in \hyperref[alg:MINRES:inexactness]{Step 7} of \cref{alg:MINRES}.
\end{condition}
\begin{condition}[$\sigma$-Limited Curvature Condition]\label{cond:LC}
	At iteration $t$ of \cref{alg:MINRES}, if $\langle\f{r}_k^{(t-1)}, \inH_k\f{r}_k^{(t-1)}\rangle \leq \sigma d\|\f{r}_k^{(t-1)}\|^2$, for some positive constant $\sigma$, then $ \mathbf{r}_k^{(t-1)} $ is declared as an $\sigma$-LC direction. This condition is effortlessly verified in \hyperref[alg:MINRES:LC]{Step 11} of \cref{alg:MINRES}. 
\end{condition}

\cref{cond:inexactness} relates the optimality of the sub-problem, i.e., it measures how close the approximated direction is to Newton's direction in terms of the residual of the normal equation. We emphasize that \cref{cond:inexactness} guarantees to be eventually satisfied during the iterations, regardless of the choice of $\eta$. This is because, as a consequence of \cite[Lemma 3.1]{liu2022minres}, $\|\inH_k\f{r}_k^{(t-1)}\|$ is decreasing to zero while $\|\inH_k\f{s}_k^{(t-1)}\|$ is monotonically increasing. This is in sharp contrast to the typical relative residual condition $ \|\f{r}_{k}^{(t-1)}\| \leq \eta \|\f{g}_{k}\| $, which is often used in related works. In nonconvex settings, where the gradient might not lie entirely in the range of Hessian, the relative residual condition might never be satisfied unless the inexactness parameter is appropriately set, which itself depends on some unavailable lower bound. 

Note that, when $\sigma = 0$,  \cref{cond:LC} coincides with the nonpositive curvature (NPC) condition. Variants of such NPC condition have been extensively used within various optimization algorithms to generate descent directions for nonconvex problems  \cite{liu2022newton,royer2020newton,xu2020second,carmon2017convex,li2023randomized}. It turns out that for the settings we consider in this paper, considering $\sigma > 0$ removes the need to introduce more structural assumptions on the Hessian matrix, while still allowing for the construction of descent directions. This is mainly due to the following simple observation that as long as \cref{cond:LC} has not been triggered, the underlying Krylov subspace contains vectors that better align with the eigenspace of Hessian corresponding to large positive eigenvalues, i.e., regions of sufficiently large positive curvature. 
\begin{lemma}\label{lemma:sigma_uniform_boundedness}
    Suppose \cref{cond:LC} has not yet been detected at iteration $t$. For any vector $\f{v} \in \krylov{t}{\inH_k}{\f{g}_k}$, we have $\lr{\f{v}, \inH_k\f{v}} \geq \sigma \|\f{v}\|^2$. 
\end{lemma}
\begin{proof}
    Let $\f{v} \in \krylov{t}{\inH_k}{\f{g}_k}$. Since, Span$\{\f{r}_k^{(0)},\f{r}_k^{(1)}, \cdots, \f{r}_k^{(t-1)}\} = \krylov{t}{\inH_k}{\f{g}_k}$, we can write $\f{v} = \sum_{i=0}^{t-1} c_i \f{r}_k^{(i)}$ for some set of scalars, $c_i$. Using the $\inH_k$-conjugacy of the residuals, it follows that 
    \begin{align*}
        \lr{\f{v}, \inH_k\f{v}} &= \sum_{i=0}^{t-1}c_i^2\lr{\f{r}_k^{(i)}, \inH_k\f{r}_k^{(i)}} \geq \sum_{i=0}^{t-1} c_i^2\sigma d \vnorm{\f{r}_k^{(i)}}^2 \geq \frac{\sigma d}{t+1} \vnorm{\sum_{i=0}^{t-1}c_i \f{r}_k^{(t-1)}}^2 \geq \sigma \vnorm{\f{v}}^{2},
    \end{align*}
    where the second-to-last inequality follows from \cite[Fact 9.7.9]{bernstein2009matrix}.
\end{proof}

We highlight that, unlike \cref{cond:inexactness}, \cref{cond:LC} does not relate to any sub-optimality criterion for the MINRES iterates. As it was shown in \cref{lemma:sigma_uniform_boundedness}, as long as \cref{cond:LC} has not been detected, the MINRES iterate $\f{s}_k^{(t-1)}$ from \cref{alg:MINRES} will be, to a great extend, aligned with the eigenspace corresponding to large positive eigenvalues and can thus be used to construct a suitable descent direction. However, \cref{cond:LC} indicates that the underlying Krylov subspace now contains directions of small or negative curvature, suggesting that the MINRES iterate $\f{s}_k^{(t-1)}$ might not be a good direction to follow. In such cases, a descent direction is constructed from the residual vector, $\f{r}_k^{(t-1)}$, rather than the iterate itself. 

Depending on which termination condition is satisfied first, the search direction for Newton-MR is constructed as follows: 
\begin{align*}
	\f{d}_{k}  = 
	\left\{\begin{array}{lll}
			\f{s}_{k}^{(t-1)},& &  \text{if \cref{cond:inexactness} is satisfied,}\\
            \\
			\f{r}_{k}^{(t-1)},& &  \text{if \cref{cond:LC} is satisfied.}
		\end{array}\right.
\end{align*} 
\noindent It turns out that MINRES enjoys a wealth of theoretical and empirical properties that position it as a preferred sub-problem solver for many Newton-type methods for nonconvex settings \cite{liu2022minres,lim2024conjugate}. Among them, one can show that the search direction $ \f{d}_{k} $ is always guaranteed to be a descent direction for $ f $. Indeed, Properties \cref{prop:descent_in_quadratic} and \cref{prop:residual_descent} below always guarantee $ \langle\f{d}_{k}, \f{g}_{k}\rangle < 0 $ when \cref{cond:inexactness,cond:LC} are met, respectively. The proofs of these properties are deferred to \cref{sec:appendix}. 
\begin{property}\label{property:MINRES}
    Consider the $t\th$ iteration of \cref{alg:MINRES} and let $0 \leq j \leq i \leq t-1$. We have
    \begin{align}
        \vnorm{\f{r}_k^{(i)}}^2 & \geq \vnorm{\f{r}_k^{(j)}}^2,\label{prop:monotonicity}\\
        -\lr{\f{r}^{(i)}_k,\f{g}_k} & = \vnorm{\f{r}^{(i)}_k}^2,\label{prop:residual_descent}\\
        \lr{\f{r}^{(t-1)}_k,\f{r}^{(i)}_k} & = \vnorm{\f{r}^{(t-1)}_k}^2.\label{prop:inner_residual}
    \end{align}
    Furthermore, if \cref{cond:LC} has not yet been detected, then 
    \begin{align}
        &\lr{\f{s}_k^{(i+1)}, \f{g}_k}  < - \lr{\f{s}_k^{(i+1)}, \inH_k\f{s}_k^{(i+1)}}\label{prop:descent_in_quadratic}\\
        &\vnorm{\f{g}_k}^2\frac{\lr{\f{g}_k, \inH\f{g}_k}^2}{\vnorm{\inH_k\f{g}_k}^4} = \vnorm{\f{s}_k^{(1)}}^2  \leq \vnorm{\f{s}_k^{(j+1)}}^2 \leq \vnorm{\f{s}_k^{(i+1)}}^2.\label{prop:increasing_s}
    \end{align}
\end{property}

Having characterized the search direction, we then enforce sufficient descent for each iteration by choosing a step size $\alpha_k$ that satisfies the Armijo-type line-search with parameter $0 < \rho < 1$ as
\begin{equation}\label{eq:armijo_linesearch}
	f(\f{x}_k + \alpha_k\f{d}_k) \leq f(\f{x}_k) + \rho\alpha_k \langle\f{g}_k, \f{s}_k\rangle.
\end{equation}
In particular, when $\f{d}_k = \f{s}_k^{(t-1)}$ (i.e., \texttt{Type = SOL}), we perform the backward tracking line-search (\cref{alg:back_tracking_ls}) to obtain the largest $0 < \alpha_k \leq 1$. Such backtracking technique is standard within the optimization methods that employ line-search as the globalization strategy \cite{nocedal2006numerical}. When $\f{d}_k = \f{r}_k^{(t-1)}$ (i.e., \texttt{Type = LC}), we perform the forward and backward tracking line-search (\cref{alg:forward_back_tracking_ls}) to find the largest $0 < \alpha_k < \infty$. After a step size, $ \alpha_{k} $ is found, the next iterate is then given as $ \f{x}_{k+1} = \f{x}_{k} + \alpha_{k} \f{d}_{k} $.
Such a forward tracking strategy, though less widely used than backward tracking, has been shown to yield much larger step sizes when directions of small or negative curvature are used and can significantly improve empirical performance \cite{gould2000exploiting,liu2022newton}. We note that even though this forward tracking strategy can also be used with $\f{d}_k = \f{s}_k^{(t-1)}$, we have found that, empirically, the larger step sizes from forward tracking when \texttt{Type = SOL} might be less effective, and at times even detrimental, compared with the case when \texttt{Type = LC}. 

The culmination of these steps is depicted in \cref{alg:NewtonMR}. Compared with the variant of Newton-MR in \cite{liu2022newton},  \cref{alg:NewtonMR} contains two main differences: employing the inexact Hessian $\inH$ instead of the exact one $\eH$, as well as leveraging $\sigma$-LC condition with $\sigma > 0$ instead of the typical NPC condition with $\sigma = 0$.  
For the sake of brevity, we have introduced a variable ``\texttt{Type}'' in \cref{alg:NewtonMR} and in our subsequent discussions. This variable takes the value ``\texttt{Type = LC}'' when \cref{cond:LC} is met, signifying $ \f{d}_{k} = \f{r}_{k}^{(t-1)} $,  and is set to ``\texttt{Type = SOL}'' when \cref{cond:inexactness} is triggered, indicating $\f{d}_{k} = \f{s}_{k}^{(t-1)}$.

\begin{algorithm}[htbp]
	\caption{Newton-MR with Inexact Hessian}
	\label{alg:NewtonMR}
	\begin{algorithmic}[1]
		\Require Initial point: $ \f{x}_0 $, First order optimality tolerance: $ 0 < \varepsilon_{\f{g}} \leq 1 $, Sub-problem inexactness tolerance : $ \eta > 0 $, Line-search parameter : $0 < \rho < 1$, Limited curvature parameter : $ \sigma > 0 $
		\State $ k = 0 $
		\While{$ \|\f{g}_k\| > \varepsilon_{\f{g}} $}
		\State Call MINRES \cref{alg:MINRES}: $(\f{d}_k, \texttt{Type}) = \text{MINRES}(\inH_k, \f{g}_k, \eta, \sigma)$ 
		\State Use \cref{alg:back_tracking_ls} (\texttt{Type = SOL}) or \cref{alg:forward_back_tracking_ls} (\texttt{Type = LC}) to obtain $\alpha_k$ satisfying \cref{eq:armijo_linesearch}
		\State $ \f{x}_{k+1} = \f{x}_k + \alpha_k \f{d}_k $
		\State $ k = k + 1 $
		\EndWhile
		\State \Return $\f{x}_k$
	\end{algorithmic}
\end{algorithm}
\begin{algorithm}[htbp]
    \caption{Backward Tracking Line-Search \cite{liu2022newton}}
    \label{alg:back_tracking_ls}
    \begin{algorithmic}[1]
        \Require $\alpha = 1$ and $0 < \xi < 1$
        \While{\cref{eq:armijo_linesearch} is not satisfied}
        \State $\alpha = \xi\alpha$
        \EndWhile
        \State \Return $\alpha$
    \end{algorithmic}
\end{algorithm}
\begin{algorithm}[htbp]
    \caption{Forward/Backward Tracking Line-Search \cite{liu2022newton}}
    \label{alg:forward_back_tracking_ls}
    \begin{algorithmic}[1]
        \Require $\alpha > 0$ and $0 < \xi < 1$
        \If{\cref{eq:armijo_linesearch} is not satisfied}
        \State Call \cref{alg:back_tracking_ls}
        \Else
        \While{\cref{eq:armijo_linesearch} is satisfied}
        \State $\alpha = \alpha/\xi$
        \EndWhile
        \State \Return $\xi\alpha$
        \EndIf
    \end{algorithmic}
\end{algorithm}

\section{Theoretical Analyses}\label{sec:theoretical_analyses}
In this section, we present complexity analyses of \cref{alg:NewtonMR}. In particular, after establishing the iteration complexity in  \cref{sec:iteration}, we derive the operation complexities in \cref{sec:operation}, which are more representative of the actual costs in the large-scale problems. \emph{Remarkably, we show that, under minimal assumptions, our algorithm converges irrespective of the degree of Hessian approximation, i.e., $\varepsilon$, and the accuracy of the inner problem solution, i.e., $\eta$ and $\sigma$.} In fact, for all of our analysis, we only make the following blanket assumption, which is widely used throughout the literature.

\begin{assumption}
	\label{assmpt:Lg}
	The function $ f $ is twice continuously differentiable and bounded below. Furthermore, there exists a constant $0 \leq \Lg < \infty$ such that for any $\f{x}, \f{y} \in \mathbb{R}^{d}$, we have $\vnorm{\f{g}(\f{x}) - \f{g}(\f{y})} \leq \Lg \vnorm{\f{x} - \f{y}}$ (or, equivalently, for any $\f{x} \in \mathbb{R}^d$, $ \vnorm{\eH(\f{x})} \leq \Lg $).
\end{assumption}

It is well-known that \cref{assmpt:Lg} implies
\begin{align}
    f(\f{x} + \f{v}) \leq f(\f{x}) + \lr{\f{g}(\f{x}), \f{v}} + \frac{\Lg}{2}\|\f{v}\|^2, \quad \forall \f{x} \in \mathbb{R}^d, \; \forall \f{v} \in \mathbb{R}^d \label{eq:second_order_taylor}.
\end{align}
Furthermore, from \cref{assmpt:Lg,eq:noise_model} and using the reverse triangle inequality, we have
\begin{align}
    \vnorm{\inH\f{v}(\f{x})} & \leq \vnorm{\inH(\f{x})}\vnorm{\f{v}} \leq (\Lg + \varepsilon)\vnorm{\f{v}}, \quad \forall \f{x} \in \mathbb{R}^d, \; \forall \f{v} \in \mathbb{R}^d. \label{eq:noisy_bound}
\end{align}

\subsection{Iteration Complexity}\label{sec:iteration}
To obtain the operation complexity of \cref{alg:NewtonMR}, we first establish its convergence rate, i.e., the iteration complexity. We emphasize that our main focus in this paper is to study the consequences of using inexact Hessian as it relates to complexity, under minimal assumptions and with minimal modifications to the original exact algorithm. Hence, we only consider minimal smoothness, as stated in \cref{assmpt:Lg}, i.e., we merely assume that the gradient is Lipschitz continuous and do not extend this smoothness assumption to the Hessian itself. Additionally, we do not employ the Hessian regularization techniques used in \cite{yao2023inexact, royer2020newton, li2023randomized}, as distorting the curvature information further could potentially lead to inferior performance; see the experiments in \cref{sec:hessiam_damp}. 

Our analysis in this section relies on \cref{lemma:linear_rate_sol,lemma:linear_rate_LC} below that shows the minimum amount of descent \cref{alg:NewtonMR} can achieve going from $ f_{k}$ to $ f_{k+1}$, when ``\texttt{Type = SOL}'' and ``\texttt{Type = LC}'', respectively.
\begin{lemma}
	\label{lemma:linear_rate_sol}
	Under \cref{assmpt:Lg}, in \cref{alg:NewtonMR} with \texttt{Type = SOL}, we have
	\begin{align*}
		f_{k+1} - f_{k} \leq -2\rho(1 - \rho) \xi C_s\|\f{g}_k\|^2,
	\end{align*}
    where $0 < \rho < 1$ and $0 < \xi < 1$ are the line-search parameters and $C_{s} \triangleq d^2 \sigma^{4} / (\Lg (\Lg + \varepsilon)^4)$.
\end{lemma}
\begin{proof}
	In this case, $\f{d}_{k} = \f{s}_{k}^{(t-1)}$. Using \cref{eq:second_order_taylor}, we have $f(\mathbf{x} + \alpha_k \f{d}_k) - f(\mathbf{x}) \leq  \alpha_k \dotprod{\f{d}_k, \f{g}_k} + \alpha_k^2 L_{\f{g}}\vnorm{\f{d}_k}^{2}/2$.
	Now, the line-search is satisfied, if, for some $ 0 < \alpha\leq 1 $, we have $\alpha_k\langle\f{d}_k,\f{g}_k\rangle + \alpha_k^2 L_{\f{g}}\|\f{d}_k\|^2/2 \leq \alpha_k\rho\langle\f{d}_k,\f{g}_k\rangle$. This inequality is satisfied for any step-size that is smaller than $- 2(1-\rho)\langle\f{d}_k,\f{g}_k\rangle / (L_{\f{g}}\|\f{d}_k\|^2)$. As a result, starting from a sufficiently large initial trial step-size, the backtracking line-search will return a step-size that is at least $\alpha_k \geq - 2(1-\rho) \xi \langle\f{d}_k,\f{g}_k\rangle / (L_{\f{g}}\|\f{d}_k\|^2)$. Since \cref{cond:LC} was not met at iteration $t-2$ (otherwise, we would have returned \texttt{Type = LC} in the previous iteration), from \cref{prop:descent_in_quadratic,lemma:sigma_uniform_boundedness}, we have $\left\langle\f{d}_k, \f{g}_k\right\rangle < - \left\langle\f{d}_k, \inH_k\f{d}_k\right\rangle \leq -\sigma \left\|\f{d}_k\right\|^2.$
	Using this and substituting back to the inequality from the line-search, we get
	\begin{align*}
		f_{k+1} - f_{k} \leq \alpha_k\rho\langle\f{d}_k,\f{g}_k\rangle \leq - 2\rho(1 - \rho)\xi \langle\f{d}_k,\f{g}_k\rangle^2/(L_{\f{g}}\|\f{d}_k\|^2) \leq - 2\sigma^2\rho(1 - \rho)\xi\|\f{d}_k\|^2/L_{\f{g}}.
	\end{align*}
    To find the lower bound for $\|\f{d}_k\|^2$, from \cref{prop:increasing_s} we get
    \begin{align*}
        \|\f{d}_k\|^2 & \geq {\lr{\f{g}_k,\inH\f{g}_k}^2}\vnorm{\f{g}_k}^2/\vnorm{\inH\f{g}_k}^4 \geq {\sigma^2 d^2} \|\f{g}_k\|^2/(\Lg + \varepsilon)^4,
    \end{align*}
    where the last inequality follows from \cref{cond:inexactness,eq:noisy_bound}. This gives the desired result.
\end{proof}

\begin{lemma}
	\label{lemma:linear_rate_LC}
	Under \cref{assmpt:Lg}, in \cref{alg:NewtonMR} with \texttt{Type = LC}, we have
	\begin{align*}
		f_{k+1} - f_{k} \leq -2\rho(1 - \rho) \xi C_l\|\f{g}_k\|^2,
	\end{align*}
    where $0 < \rho < 1$  and $0 < \xi < 1$ are the line-search parameters and $C_l \triangleq \eta^2/[L_{\f{g}}((\Lg + \varepsilon)^2 + \eta^2)]$.
\end{lemma}
\begin{proof}
	First, we note that when \texttt{Type = LC}, \cref{cond:inexactness} has not been met (since it appears before \cref{cond:LC} in \cref{alg:MINRES}). Similar to the proof \cref{lemma:linear_rate_sol}, for $ \alpha_{k} $ to satisfy \cref{eq:armijo_linesearch}, we need $\alpha_k\langle\f{d}_k,\f{g}_k\rangle + \alpha_k^2L_{\f{g}}\|\f{d}_k\|^2 / 2 \leq \alpha_k\rho\langle\f{d}_k,\f{g}_k\rangle$, which implies any step-size smaller than $-2(1-\rho)\langle\f{d}_k,\f{g}_k\rangle / (L_{\f{g}}\|\f{d}_k\|^2)$ satisfies the line search. Since \cref{cond:LC} is satisfied, i.e., $\f{d}_k = \f{r}_k^{(t-1)}$, from \cref{prop:residual_descent}, we have $\langle\f{d}_k ,\f{g}_k\rangle = \langle\f{r}_k^{(t-1)},\f{g}_k\rangle = - \|\f{r}_k^{(t-1)}\|^2 = -\vnorm{\f{d}_k}^2$.  Hence, using a similar argument as in the proof of \cref{lemma:linear_rate_sol}, the step size returned from the line-search must satisfy $\alpha_k \geq 2(1 - \rho) \xi / L_{\f{g}}$. Substituting back to the inequality from the line-search, we get $f_{k+1} - f_k \leq \alpha_k\rho\langle\f{d}_k,\f{g}_k\rangle = -2\rho(1 - \rho )\xi\vnorm{\f{d}_k}^2 / L_{\f{g}}$.
	Since \cref{cond:inexactness} is not satisfied, we have
	\begin{align*}
		\frac{\vnorm{\inH_k\f{r}_k^{(t-1)}}^{2}}{\eta^{2}} & \geq \vnorm{\inH_{k} \f{s}_{k}^{(t-1)}}^{2} = \vnorm{\f{r}_k^{(t-1)} + \f{g}_{k}}^{2} \\
  &= \vnorm{\f{r}_k^{(t-1)}}^2 + \vnorm{\f{g}_k}^2 + 2\lr{\f{r}_k^{(t-1)},\f{g}_k} = \vnorm{\f{g}_{k}}^{2} - \vnorm{\f{r}_k^{(t-1)}}^{2},
	\end{align*}
	where the last equality follows from \cref{prop:residual_descent}. Together, this implies $((L_{\f{g}} + \varepsilon)^{2} + \eta^{2})\|\f{r}_{k}^{(t-1)}\|^{2} \geq \|\inH_k\f{r}_k^{(t-1)}\|^2 + \eta^2 \|\f{r}_k^{(t-1)}\|^2 \geq \eta^2\vnorm{\f{g}_k}^2$,
	which gives the desired result.
\end{proof}

With \cref{lemma:linear_rate_sol,lemma:linear_rate_LC} in hand, we can study the convergence rate of \cref{alg:NewtonMR} for nonconvex functions.  Beyond general nonconvex settings, we also consider subclass of functions that satisfy the $\theta$-PL condition, i.e., \cref{def:PL}. In this setting, we show that with $\theta = 2$ and under \cref{assmpt:Lg}, \cref{alg:NewtonMR} achieves a global linear convergence rate. For  $1 \leq \theta < 2$, however, the global convergence of \cref{alg:NewtonMR} transitions from a linear rate to a sub-linear rate in the neighbourhood of a solution. 
In all these cases, $\theta$-PL property allows us to derive a faster convergence rate than in the general nonconvex settings.  
In the context of Newton-type methods with inexact Hessian information, to the best of our knowledge, our work here is the first to establish fast convergence rates for $ \theta $-PL functions under mild smoothness assumptions and algorithmic requirements.

Prior works on Newton-type methods with inexact Hessians and linear convergence rates have predominantly focused on (strongly) convex functions. These studies have either assumed that the approximate Hessian remains positive definite \cite{roosta2019sub,byrd2011use,byrd2012sample,bollapragada2016exact,erdogdu2015convergence,roosta2019sub} or that the gradient lies in the range of Hessian at all times \cite{gower2019rsn}. These stringent assumptions either necessitate highly accurate Hessian approximations or significantly restrict the scope of problems that can be effectively solved.
However, these limitations often conflict with what is observed in practice. For instance, the sub-sampled Newton method can effectively solve (strongly) convex finite-sum problems with very crude Hessian approximations using much fewer samples than what the concentration inequalities predict \cite{xu2016sub}. Similarly, for high-dimensional convex problems where the Hessian matrix is singular and the gradient may not entirely lie in the range of Hessian, Newton's method remains practically a viable algorithm \cite{roosta2022newton}.

The limitations imposed by these stringent assumptions arise from the common practice in previous studies, where the sub-problem solver is treated as a black box, neglecting the opportunity to leverage its inherent properties.
%
%
In contrast, by tapping into the properties of MINRES, we show that simply requiring \cref{assmpt:Lg} allows us to establish the convergence of \cref{alg:NewtonMR} for nonconvex functions under \emph{any amount of inexactness} in either the Hessian matrix or the sub-problem solution, effectively bridging the gap between theory and practice in these setting.

\begin{theorem}
    \label{thm:pl}
    Let \cref{assmpt:Lg} holds and consider \cref{alg:NewtonMR} with any $ \varepsilon > 0 $ in \cref{eq:noise_model}, any $\eta > 0$ in \cref{cond:inexactness}, any $\sigma > 0$ in \cref{cond:LC}, the line-search parameter $0 < \rho < 1$. 
    \begin{itemize}
        \item For functions satisfying the $\theta$-PL condition \cref{cond:PL}, after at most
    \begin{align*}
        K \triangleq \frac{1}{2\mu C} \left(\log \left(\max\{1,f_0 - f^\star\}\right) + \varepsilon_{f}^{(1-{2}/{\theta})} \log (1/\varepsilon_{f})\right),
    \end{align*}
    iterations of Algorithm \ref{alg:NewtonMR}, the approximate global optimality \cref{eq:global_optimality_condition} is satisfied. 
    \item For more general nonconvex functions, after at most
    \begin{align*}
        K \triangleq \frac{(f_0 - f^\star)\varepsilon_{\f{g}}^{-2}}{C},
    \end{align*}
    iterations of Algorithm \ref{alg:NewtonMR}, the approximate first order sub-optimality \cref{eq:1st_order_optimality_condition} is satisfied.
    \end{itemize}
    Here, $C \triangleq 2\rho(1 - \rho) \xi \min\left\{C_s, C_l\right\}$ where $C_s$ and $C_l$ are defined in \cref{lemma:linear_rate_sol,lemma:linear_rate_LC}, respectively.
\end{theorem}
\begin{proof}
    For each iteration of \cref{alg:NewtonMR}, we either have \texttt{Type = SOL} or \texttt{Type = LC}. 
    Hence, from \cref{lemma:linear_rate_sol,lemma:linear_rate_LC}, we obtain
    \begin{align}
		f_{k+1} - f_k \leq -C\|\f{g}_k\|^2.\label{eq:min_reduction}
	\end{align}
    Suppose the $\theta$-PL condition \cref{cond:PL} holds. When $\theta = 2$, it immediately follows that
    \begin{align*}
        f_{k+1} - f^\star \leq \left(1 - 2\mu C\right)(f_{k} - f^\star).
    \end{align*}
    Hence, to get $f_{k} - f^\star \leq \varepsilon_{f}$, it suffices to find $k$ for which $\left(1 - 2\mu C\right)^{k} \leq \varepsilon_{f}/(f_{0}-f^\star)$, which gives the desired results. Now, consider $1 \leq \theta < 2$. As long as $f_k - f^\star \geq 1$, we have $\vnorm{\f{g}_k}^2/2\mu \geq (f_{k} - f^\star)^{2/\theta} \geq f_{k} - f^\star$, which is handled just as in the case of $\theta = 2$, implying that after at most $\log (f_0 - f^\star)/(2\mu C)$ iterations, we have $f_k - f^\star \leq 1$. Once, $\varepsilon_{f} \leq f_k - f^\star < 1$, from \cref{eq:min_reduction,cond:PL}, we obtain
    \begin{align*}
        f_{k+1} - f^\star & \leq \left(1 - \frac{C \|\f{g}_k\|^2}{f_k - f^\star}\right)(f_k - f^\star) \leq \left(1 - 2\mu C(f_k - f^\star)^{{2}/{\theta} - 1}\right)(f_k - f^\star) \\
        & \leq \left(1 - 2\mu C\varepsilon_{f}^{\frac{2}{\theta} - 1}\right)(f_k - f^\star).
    \end{align*}
    Setting $ (1 - 2\mu C\varepsilon_{f}^{\frac{2}{\theta} - 1})^{k} \leq \varepsilon_{f} $ gives the result.

    
    Now, consider a more general nonconvex setting. By a way of contradiction, supposed $\vnorm{\f{g}_{J+1}} \leq \varepsilon_{\f{g}}$ happens for the first time at some $J > K$. So, we must have $\|\f{g}_k\| > \varepsilon_{\f{g}}$, for all $j \leq J$. From \cref{eq:min_reduction}, it follows that 
    \begin{align*}
        f_{j} - f_{j+1} & \geq C \|\f{g}_j\|^2 > C \varepsilon_{\f{g}}^2, \quad \forall j \leq J. 
    \end{align*}
    Hence, by the definition of $K$, we must have
    \begin{align*}
        f_0 - f_{J} & = \sum_{j=0}^{J-1} f_j - f_{j+1} > J C \varepsilon_{\f{g}}^2  > K C \varepsilon_{\f{g}}^2 = f_0 - f^\star
    \end{align*}
    This implies $f^\star > f_{J}$, which is a contradiction. 
\end{proof}

\begin{remark}
    As it can be seen, for $\theta = 2$, the convergence is linear, while for $1\leq \theta < 2$, the convergence transitions to a sublinear rate in the vicinity of an optimum point, i.e., when $f_{k} - f^{\star} \leq 1$. Smaller values of $\theta$ imply a slower convergence rate in the neighborhood of a solution, which is expected as the function becomes flatter in those regions. 
\end{remark}

\begin{remark}
The advantage of \cref{thm:pl} is that it establishes linear convergence of \cref{alg:NewtonMR} for \emph{any} value of $ \varepsilon $ and $ \eta $. In other words, regardless of how crude our Hessian approximation or inner problem solution might be, \cref{alg:NewtonMR} still converges for nonconvex functions. This is only possible thanks to leveraging the properties of MINRES as a sub-problem solver in \cref{lemma:linear_rate_LC,lemma:linear_rate_sol}. This, to a great degree, corroborates what is typically observed in practice. For example, for solving strongly convex problems, which can be seen as a sub-class of PL functions, Newton's method remains a convergent algorithm with well-defined iterations, irrespective of the level of approximation used for Hessian or the underlying sub-problem solution \cite{roosta2022newton,roosta2019sub,berahas2020investigation}. 
\end{remark}	

\begin{remark}
    The iteration complexity of $\mathcal{O}(\varepsilon_{\f{g}}^{-2})$ in \cref{thm:pl} is known to be optimal within the class of second-order methods when applied to solve \cref{eq:min_f} under Lipschitz continuous gradient assumption \cite{cartis2018worst}. Therefore, under \cref{assmpt:Lg}, our analysis aligns with the theoretical lower bound, highlighting the efficiency of our approach under these conditions. 
\end{remark}

\subsection{Operation Complexity}\label{sec:operation}
In large-scale settings where storing the Hessian matrix explicitly is impractical, and one can only employ operations such as matrix-vector products, the main computational cost lies in solving the sub-problem. Therefore, iteration complexity alone does not adequately represent the actual costs. In this section, we provide an upper bound estimate on the operation complexity, that is the total number of gradient and Hessian-vector product evaluations required by \cref{alg:NewtonMR}, to achieve a desired sub-optimality.
%
For clarity of exposition, we drop the subscript $k$ of $\f{g}_k$ and $\inH_k$ in what follows in this section.

The interplay between the inexact Hessian and the gradient plays a central role in establishing the operation complexity of \cref{alg:NewtonMR}. This interplay is captured by the notion of $\f{g}$-relevant eigenvalues/eigenvectors.
\begin{definition}
    An eigenvector $\f{v}$ of $\inH$ is called a $\f{g}$-relevant eigenvector if $\lr{\f{g}, \f{v}} \not= 0$. The corresponding eigenvalue $\lambda$, i.e., $\inH\f{v} = \lambda\f{v}$, is called a $\f{g}$-relevant eigenvalue.
\end{definition}

Consider the eigen-decomposition of $\inH$ be
\begin{align}
    \label{eq:eigen_decomp}
    \inH = \begin{bmatrix}
        \f{V} & \f{V}_\perp & \f{V}_\f{0}
        \end{bmatrix}\begin{bmatrix}
        \f{Z} & & \\
        & \f{Z}_\perp & \\
        & & \f{0}
    \end{bmatrix}\begin{bmatrix}
        \f{V} & \f{V}_\perp & \f{V}_\f{0}
        \end{bmatrix}^\top,
\end{align}
where $\f{Z} = \text{diag}(\lambda_1, \cdots, \lambda_{\phi})$ is the diagonal matrix with \textit{nonzero} $\f{g}$-relevant eigenvalues in a nonincreasing order, that is $\lambda_1 \geq \cdots \geq \lambda_{\phi_+} > 0 > \lambda_{\phi_{+}+1} \geq \cdots \geq \lambda_{\phi}$, where $\phi$ and $\phi_+$ denote, respectively, the total number of nonzero and positive $\f{g}$-relevant eigenvalues (including multiplicities). The columns of the matrix $\f{V}$ contain the corresponding orthonormal eigenvectors. The matrices $\f{V}_\perp$ and $\f{Z}_\perp$ are defined similarly, but with respect to \textit{nonzero} non-$\f{g}$-relevant counterparts, i.e., $\f{V}_\perp \f{g} = \f{0}$. Finally, the columns of the matrix $\f{V}_\f{0}$ span the nullspace of $\inH$. 
For any $1 \leq j \leq \phi$, we further define $\f{V} = \begin{bmatrix} \f{V}_j & \f{V}_{j\perp} \end{bmatrix}$, where $\f{V}_j$ denotes the matrix formed by the first $j$ columns of $\f{V}$ and $\f{V}_{j\perp}$ contains the remaining columns of $\f{V}$.
Similarly, we define $\f{Z}_j = \text{diag}(\lambda_1, \cdots, \lambda_{j})$ and $\f{Z}_{j\perp} = \text{diag}(\lambda_{j+1}, \cdots, \lambda_{\phi})$. A simple example can help clarify the notation.
\begin{example}
    Let $\f{U} = \begin{bmatrix}
        \f{u}_1 & \cdots & \f{u}_{10}
    \end{bmatrix} \in \mathbb{R}^{10\times 10}$ be an orthonormal matrix, $\f{g} = \f{u}_1 + \f{u}_3 + \f{u}_7 + \f{u}_{10}$ and
    \begin{align*}
        \inH & = \f{U} \; \text{diag}(9,7,6,6,4,4,0,-2,-3,-5)\;\f{U}^\top.
    \end{align*}
    Then, $\f{V} = \begin{bmatrix}
        \f{u}_1 & \f{u}_3 & \f{u}_4 & \f{u}_{10}
    \end{bmatrix}$, $\f{V}_\perp = \begin{bmatrix}
        \f{u}_2 & \f{u}_5 & \f{u}_6 & \f{u}_{8} & \f{u}_{9}
    \end{bmatrix}$, $\f{V}_0 = \f{u}_7$; $\lambda_1 = 9$, $\lambda_2 = 3$, $\lambda_3 = 3$ and $\lambda_4 = -5$; $\phi_+ = 3$ and $\phi = 4$.
\end{example}

We now make several key observations when \cref{cond:LC} is not met in the very first iteration. In this case, from \cref{eq:noisy_bound}, we must have $\sigma d < (\Lg + \varepsilon) $. Also, clearly, $\phi_+ \geq 1$. Furthermore, we must have have $\lambda_j > \sigma d / 2$, for some $1 \leq j \leq \phi_+$. 
Indeed, suppose the contrary, which implies $\lambda_{1} \leq \sigma d / 2$. 
Since \cref{cond:LC} is not detected in the first iteration, i.e., $\lr{\f{g}, \inH\f{g}} > \sigma d \vnorm{\f{g}}^2$, we get
\begin{align*}
    \sigma d < \frac{\lr{\f{g},\inH\f{g}}}{\vnorm{\f{g}}^2} = \sum_{i=1}^\phi\lambda_i \frac{\lr{\f{g},\f{v}_i}^2}{\|\f{g}\|^2} \leq \lambda_{1} \sum_{i=1}^\phi \frac{\lr{\f{g},\f{v}_i}^2}{\vnorm{\f{g}}^2} \leq \lambda_{1} \leq \frac{\sigma d}{2},
\end{align*}
which is a contradiction. Now, let $1 \leq j \leq \phi_+$ be the largest index such that $\lambda_j > \sigma d / 2$, i.e., $\lambda_j$ is the smallest $\f{g}$-relevant eigenvalue with such property. Since $\sigma d\|\f{g}\|^2 < \lr{\f{g}, \inH\f{g}}$, the Pythagorean theorem gives 
\begin{align*}
    \sigma d\vnorm{(\f{VV}^\top - \f{V}_{j}\f{V}_{j}^\top)\f{g}}^2 & + \sigma d \vnorm{\f{V}_{j}\f{V}_{j}^\top\f{g}}^2 + \sigma d\vnorm{\f{V}_\perp\f{V}_\perp^\top\f{g}}^2\\
    & < \lr{\f{V}_{j}\f{V}_{j}^\top\f{g}, \inH\f{V}_{j}\f{V}_{j}^\top\f{g}} + \lr{(\f{VV}^\top - \f{V}_{j}\f{V}_{j}^\top)\f{g}, \inH(\f{VV}^\top - \f{V}_{j}\f{V}_{j}^\top)\f{g}} \\
    & \leq (\Lg + \varepsilon) \vnorm{\f{V}_{j}\f{V}_{j}^\top\f{g}}^2 + \frac{\sigma d}{2}\vnorm{(\f{VV}^\top - \f{V}_{j}\f{V}_{j}^\top)\f{g}}^2,
\end{align*}
where the last inequality follows from \cref{eq:noisy_bound} and the fact that $\f{V}_{j\perp} \f{V}_{j\perp}^{\top} = \f{VV}^\top - \f{V}_{j}\f{V}_{j}^\top$ is the projection on the eigenspace corresponding to the eigenvalues that are smaller than or equal to $\sigma d / 2$. Hence, it follows that 
\begin{align}\label{eq:eigenvector_bound}
    \frac{\sigma d}{2}\vnorm{(\f{VV}^\top - \f{V}_{j}\f{V}_{j}^\top)\f{g}}^2 & \leq (\Lg + \varepsilon - \sigma d)\vnorm{\f{V}_{j}\f{V}_{j}^\top\f{g}}^2.
\end{align}

Next, we state a sufficient condition for \cref{cond:inexactness} to be satisfied.
\begin{lemma}\label{lemma:sufficient_inexactness}
    Suppose \cref{assmpt:Lg} holds and $\phi_+ \geq 1$. For any $1 \leq j \leq \phi_+$, if
    \begin{align*}
        \vnorm{\f{r}^{(t-1)}}^2 & \leq \frac{\eta^2}{(\Lg + \varepsilon)^2 + \eta^2}\vnorm{\f{V}_j\f{V}_j^\top\f{g}}^2,
    \end{align*}
    then \cref{cond:inexactness} is satisfied.
\end{lemma}
\begin{proof}
    By assumption on $\|\f{r}^{(t-1)}\|$, We have
    \begin{align*}
        \vnorm{\f{r}^{(t-1)}}^2 & \leq \frac{\eta^2}{(\Lg + \varepsilon)^2 + \eta^2}\vnorm{\f{V}_j\f{V}_j^\top\f{g}}^2 \leq \frac{\eta^2}{(\Lg + \varepsilon)^2 + \eta^2}\vnorm{\f{g}}^2,
    \end{align*}
    which using  \cref{eq:noisy_bound} implies
    \begin{align*}
        \vnorm{\inH\f{r}^{(t-1)}}^2 + \eta^2\vnorm{\f{r}^{(t-1)}}^2 & \leq \left((\Lg + \varepsilon)^2 + \eta^2\right)\vnorm{\f{r}^{(t-1)}}^2 \leq \eta^2\vnorm{\f{g}}^2.
    \end{align*}
    Hence, using \cref{prop:residual_descent} and the fact that $\|\f{g}\|^{2} - \|\f{r}^{(t-1)}\|^{2} = \|\f{g}\|^{2} + 2\langle\f{g},\f{r}^{(t-1)}\rangle + \|\f{r}^{(t-1)}\|^{2} = \|\f{g} + \f{r}^{(t-1)}\|^2$, we obtain
    \begin{align*}
        \vnorm{\inH\f{r}^{(t-1)}}^2 & \leq \eta^2\left( \vnorm{\f{g}}^2 - \vnorm{\f{r}^{(t-1)}}^2 \right) = \eta^2\vnorm{\f{g} + \f{r}^{(t-1)}}^2 = \eta^2\vnorm{\inH\f{s}^{(t-1)}}^2.
    \end{align*}
\end{proof}
The following lemma gives a certain convergence behaviour of MINRES for indefinite matrices, as it relates to a subspace corresponding to positive $\f{g}$-relevant eigenvalues. 
\begin{lemma}\label{lemma:minres_convergence}
    Let $\phi_+ \geq 1$ and $1 \leq i \leq \phi_+$. For any $0 < \omega < 4$, after at most 
    \begin{align*}
        t \geq \left\lceil\frac{\sqrt{\lambda_1/\lambda_i}}{4}\log\left(\frac{4}{\omega}\right)\right\rceil + 1,
    \end{align*}
    iterations of  \cref{alg:MINRES}, we have $\|\f{r}^{(t-1)}\|^2 \leq \|(\f{VV} - \f{V}_{i}\f{V}_{i}^\top)\f{g}\|^2 + \omega\|\f{V}_{i}\f{V}_{i}^\top\f{g}\|^2$.
\end{lemma}
\begin{proof}
Recall that the $ t\th $ iterate of MINRES using inexact Hessian $\inH(\f{x})$ is given as \cite{liu2022minres,paige1975solution} 
\begin{align}\label{eq:MINRES_min_problem}
	\f{s}_k^{(t)}  = \argmin_{\f{s} \in \mathcal{K}_{t}(\inH_k, \f{g}_k)}  \vnorm{\inH_k \f{s} + \f{g}_k}^{2}.
\end{align}
From \cref{eq:eigen_decomp}, $\text{Range}(\inH) = \text{Range}(\f{V}_{i}\f{V}_{i}^\top) \bigoplus \text{Range}(\f{V}_{i\perp}\f{V}_{i\perp}^\top + \f{V}_{\perp}\f{V}_{\perp}^\top)$, for any $1 \leq i \leq \phi_+$. Using a similar reasoning as \cite[Lemma 3.12]{liu2021convergence},  the formulation \cref{eq:MINRES_min_problem} can be written as
\begin{align*}
    \vnorm{\f{r}^{(t)}}^2 & = \min_{\f{s}\in\mathcal{K}_t(\inH,\f{g})} \vnorm{\inH\f{s} + \f{g}}^2\\
    & = \min_{\f{s}\in\f{V}_{i}\f{V}_{i}^\top\mathcal{K}_t(\inH,\f{g})} \vnorm{\inH\f{s} + \f{V}_{i}\f{V}_{i}^\top\f{g}}^2 + \min_{\f{s}\in(\f{V}_{i\perp}\f{V}_{i\perp}^\top + \f{V}_{\perp}\f{V}_{\perp}^\top)\mathcal{K}_t(\inH,\f{g})} \vnorm{\inH\f{s} + (\f{V}_{i\perp}\f{V}_{i\perp}^\top + \f{V}_{\perp}\f{V}_{\perp}^\top)\f{g}}^2\\
    & = \min_{\f{s}\in\f{V}_{i}\f{V}_{i}^\top\mathcal{K}_t(\inH,\f{g})} \vnorm{\inH\f{s} + \f{V}_{i}\f{V}_{i}^\top\f{g}}^2 + \min_{\f{s}\in\f{V}_{i\perp}\f{V}_{i\perp}^\top\mathcal{K}_t(\inH,\f{g})} \vnorm{\inH\f{s} + \f{V}_{i\perp}\f{V}_{i\perp}^\top\f{g}}^2,
\end{align*}
where the last equality follows from the fact that $\f{V}_{\perp}$ are constructed from non-$\f{g}$-relevant eigenvectors and hence, $\f{V}_{\perp}\f{V}_{\perp}^\top\f{g} = \f{0}$ and $\f{V}_{\perp}\f{V}_{\perp}^\top\mathcal{K}_t(\inH,\f{g}) = \{\f{0}\}$. We also note that
\begin{align*}
    \min_{\f{s}\in\f{V}_{i\perp}\f{V}_{i\perp}^\top\mathcal{K}_t(\inH,\f{g})} \vnorm{\inH\f{s} + \f{V}_{i\perp}\f{V}_{i\perp}^\top\f{g}}^2 \leq \vnorm{\f{V}_{i\perp}\f{V}_{i\perp}^\top\f{g}}^2 = \vnorm{(\f{VV}^\top - \f{V}_{i}\f{V}_{i}^\top)\f{g}}^2.
\end{align*}
Since for all $\f{s}\in\f{V}_{i}\f{V}_{i}^\top\mathcal{K}_t(\inH,\f{g})$, we can write $\f{s} = \f{V}_{i}\f{V}_{i}^\top\f{s}$, it follows that
\begin{align*}
    & \min_{\f{s}\in\f{V}_{i}\f{V}_{i}^\top\mathcal{K}_t(\inH,\f{g})} \vnorm{\inH\f{V}_{i}\f{V}_{i}^\top\f{s} + \f{V}_{i}\f{V}_{i}^\top\f{g}}^2 = \min_{\f{s}\in\mathcal{K}_t(\f{V}_{i}\f{Z}_i\f{V}_{i}^\top,\f{V}_i\f{V}_{i}^\top\f{g})} \vnorm{\f{V}_{i}\f{Z}_{i}\f{V}_{i}^\top\f{s} + \f{V}_{i}\f{V}_{i}^\top\f{g}}^2.
\end{align*}
where $\f{Z}_i$ is the diagonal matrix containing the top $i$ positive $\f{g}$-relevant eigenvalues of $\inH$. Hence, using the standard convergence results for MINRES \cite{saad2003iterative}, we have 
\begin{align*}
    \min_{\f{s}\in\mathcal{K}_t(\f{V}_{i}\f{Z}_i\f{V}_{i}^\top,\f{V}_i\f{V}_{i}^\top\f{g})} \vnorm{\f{V}_{i}\f{Z}_{i}\f{V}_{i}^\top\f{s} + \f{V}_{i}\f{V}_{i}^\top\f{g}}^2 & \leq 4\vnorm{\f{V}_{i}\f{V}_{i}^\top\f{g}}^2\left(\frac{\sqrt{\lambda_1 / \lambda_i} - 1}{\sqrt{\lambda_1 / \lambda_i} + 1}\right)^{2t}.
\end{align*}
Together, we obtain,
\begin{align*}
    \vnorm{\f{r}^{(t)}}^2 & \leq \vnorm{(\f{VV} - \f{V}_{i}\f{V}_{i}^\top)\f{g}}^2 + 4\vnorm{\f{V}_{i}\f{V}_{i}^\top\f{g}}^2\left(\frac{\sqrt{\lambda_1 / \lambda_i} - 1}{\sqrt{\lambda_1 / \lambda_i} + 1}\right)^{2t}.
\end{align*}
Now, for any $0 < \omega < 4$, if \cref{alg:MINRES} iterates for
\begin{align*}
    t \geq \left\lceil\frac{\sqrt{\lambda_1 / \lambda_i}}{4}\log\left(\frac{4}{\omega}\right)\right\rceil \geq \left\lceil\log\left(\frac{4}{\omega}\right)\bigg/\log\left(\frac{\sqrt{\lambda_1 / \lambda_i} - 1}{\sqrt{\lambda_1 / \lambda_i} + 1}\right)^2\right\rceil,
\end{align*}
then we have $\|\f{r}^{(t)}\|^2 \leq \vnorm{(\f{VV}^\top - \f{V}_{i}\f{V}_{i}^\top)\f{g}}^2 + \omega\vnorm{\f{V}_{i}\f{V}_{i}^\top\f{g}}^2$.
\end{proof}

Using \cref{lemma:sufficient_inexactness,lemma:minres_convergence}, the following lemma gives an upper estimate on the total number of iterations of \cref{alg:MINRES} before termination.
\begin{lemma}\label{lemma:operation_complexity}
    Suppose \cref{assmpt:Lg} holds and let $\sigma > 0$ and $\eta > 0$. \cref{alg:MINRES} terminates after at most $T$ iterations where 
\begin{align}\label{eq:minres_iteration}
    T \triangleq \left\{\begin{array}{ll}
       1,  &  \displaystyle \frac{\sigma d}{(\Lg+\varepsilon)} \geq 1,\\\\
       \displaystyle \left\lceil\sqrt{\frac{(\Lg + \varepsilon)}{8 \sigma d}} \log\left(\frac{8((\Lg + \varepsilon)^2 + \eta^2)}{\eta^2}\right)\right\rceil+1,  & \displaystyle  \frac{\eta^2 + (\Lg + \varepsilon)^2}{ 1.25 \eta^2 + (\Lg + \varepsilon)^2} \leq \frac{\sigma d}{(\Lg+\varepsilon)} < 1, \\\\
       \displaystyle \left\lceil\frac{(\Lg + \varepsilon) ((\Lg + \varepsilon)^2 + \eta^2)}{\sigma (1.25 \eta^2 + (\Lg + \varepsilon)^2)}\right\rceil, & \text{Otherwise.}
    \end{array}\right.
\end{align}
\end{lemma}
\begin{proof}
    If \cref{cond:LC} is met in the very first iteration, then the result trivially holds. Otherwise, by the observations made earlier, we must have $\sigma d < \Lg + \varepsilon$ and $\lambda_j > \sigma d / 2$, for some $1 \leq j \leq \phi_+$. Also, from \cref{eq:eigenvector_bound}, for this $j$ we get,
    $\|(\f{VV}^\top - \f{V}_{j}\f{V}_{j}^\top)\f{g}\|^2  \leq 2 (\Lg + \varepsilon - \sigma d)/(\sigma d)\|\f{V}_{j}\f{V}_{j}^\top\f{g}\|^2$. In this case, for the sufficient condition from \cref{lemma:sufficient_inexactness} to be satisfied, we need to find $\omega$ in \cref{lemma:minres_convergence} such that
    \begin{align*}
        \|(\f{VV}^\top - \f{V}_{j}\f{V}_{j}^\top)\f{g}\|^2 + \omega\|\f{V}_{j}\f{V}_{j}^\top\f{g}\|^2 \leq \frac{\eta^2}{(\Lg + \varepsilon)^2 + \eta^2}\vnorm{\f{V}_j\f{V}_j^\top\f{g}}^2.
    \end{align*}
    Hence, it suffices to find $\omega$ such that 
    \begin{align*}
        \frac{2 (\Lg + \varepsilon - \sigma d)}{\sigma d} + \omega \leq \frac{\eta^2}{(\Lg + \varepsilon)^2 + \eta^2}.
    \end{align*}
    It follows that, as long as
    \begin{align}
        \label{eq:d}
        d \geq \frac{2 (\Lg + \varepsilon) ((\Lg + \varepsilon)^2 + \eta^2)}{\sigma (2.5 \eta^2 + 2(\Lg + \varepsilon)^2)},
    \end{align}
    we have 
    \begin{align*}
        \frac{\eta^2}{(\Lg + \varepsilon)^2 + \eta^2} - \frac{2 (\Lg + \varepsilon - \sigma d)}{\sigma d} > \frac{0.5\eta^2}{(\Lg + \varepsilon)^2 + \eta^2}, 
    \end{align*}
    which implies that
    \begin{align*}
        \omega \leq \frac{0.5\eta^2}{(\Lg + \varepsilon)^2 + \eta^2},
    \end{align*}
    ensures \cref{cond:inexactness}. Otherwise, \cref{alg:MINRES} naturally terminates after at most $d $ iterations with an exact solution. Noting that $\lambda_{1} \leq \Lg + \varepsilon$ and $\lambda_j > \sigma d / 2$, \cref{lemma:minres_convergence} implies the result.
\end{proof}

With \cref{lemma:operation_complexity}, it is straightforward to obtain the total operation complexity of \cref{alg:NewtonMR}, which is again optimal among the class of second-order methods for solving non-convex problems of the form \cref{eq:min_f} satisfying Lipschitz continuous gradient assumption \cite{cartis2018worst}. 

\begin{corollary}[Optimal Operation Complexity of \cref{alg:NewtonMR}]
\label{thm:1st_ operation}
Under the assumptions of \cref{thm:pl}, for $\theta$-PL functions, after at most $\mathcal{O}\left(\log \left(\max\{1,f_0 - f^\star\}\right) + \varepsilon_{f}^{(1-{2}/{\theta})} \log (1/\varepsilon_{f})\right)$ gradient and Hessian-vector product evaluations, \cref{alg:NewtonMR} produces a solution that satisfies the approximate global optimality \cref{eq:global_optimality_condition}. For more general nonconvex functions, this bound is $\mathcal{O} \left(\varepsilon_{\f{g}}^{-2}\right)$, for the approximate first order sub-optimality \cref{eq:1st_order_optimality_condition} to be satisfied.
\end{corollary} 

\section{Numerical Experiments}\label{sec:numerical_exp}
For our numerical experiments, we consider a series of finite-sum minimization problems. We evaluate the performance of \cref{alg:NewtonMR} in comparison with several other second order methods, namely Newton-CG \cite{royer2020newton} and its sub-sampled variant \cite{yao2023inexact}, Steihaug's trust-region method \cite{conn2000trust, nocedal2006numerical} and its corresponding sub-sampled algorithms  \cite{xu2020newton,yao2021inexact}, as well as L-BFGS \cite{nocedal2006numerical}.
We approximate the Hessian by means of sub-sampling, that is we form $\inH$ using $1\%$, $5\%$, and $10\%$ of the total samples. For each experiment, the performance is depicted in two plots, namely the graph of $f(\f{x})$ with respect to the number of iterations as well as a corresponding plot based on the number of function oracle calls. For instance, the amount of work required for computing the gradient can be considered as equivalent to two oracle calls, i.e., one forward and one backward propagation, and a Hessian-vector product operation can be obtained through a work equivalent to four oracle calls (see \cite{liu2021convergence,roosta2022newton} for more details). Our implementations are matrix-free in that we do not store the  (inexact) Hessian and rely on Hessian-vector operations to perform the iterations of MINRES and CG. The code for the following experiments can be found in \url{https://github.com/alexlim1993/NumOpt}.

\subsection{Implementation details} 
\paragraph{Newton-MR and its sub-sampled variants}: For Armijo line-search parameters in \cref{eq:armijo_linesearch}, we use $\rho = 10^{-4}$, which is typical \cite{nocedal2006numerical}. For both \cref{alg:back_tracking_ls,alg:forward_back_tracking_ls}, we use $\xi = 0.5$. 
The parameters $\sigma$ in \cref{cond:LC} and $\eta$ in \cref{cond:inexactness} will be chosen specific to each problem. In particular, we aim to fine-tune $\eta$ so as to enable generating directions of \texttt{Type = SOL} more frequently.

\paragraph{Newton-CG and its sub-sampled variants}: We implement the line-search used in \cite{royer2020newton,yao2023inexact} with its parameter set to $10^{-4}$. The back-tracking line-search parameter is set to $0.5$. Within the Capped-CG algorithm \cite{royer2020newton}, the damping parameters $\tilde{\f{H}} = \eH + \zeta\f{I}$ is set to $\zeta = 0.5$. The desired accuracy of Capped-CG algorithm varies across different experiments to allow for best performance.

\paragraph{Steihaug's trust-region and its sub-sampled variants}: Following the notation used in \cite[Algorithm 4.1]{nocedal2006numerical}, we set $\hat{\Delta} = 10^{10}$, $\Delta_0 = 10^{5}$ and $\eta = 0.1$. We use the CG–Steihaug \cite[Algorithm 7.2]{nocedal2006numerical} to approximately solve the trust-region sub-problems.

\paragraph{L-BFGS}: The limited memory size is set to 20, and the Armijo and the Strong Wolfe line-search parameters are, respectively, set to their typical values of $10^{-4}$ and $0.9$.\\

All algorithms are terminated when either the oracle calls reach $10^{6}$, the magnitude of the gradient falls below $10^{-6}$, or the step size/trust region goes below $10^{-18}$. In the legends of \cref{fig:MNIST_logloss,fig:MNIST_NLS}, `$x \%$' refers to when the approximate Hessian is formed using an `$x$' percent of the full dataset. For all other figures, `\texttt{algorithm}\_$x$' refers to a given algorithm when $100 x$ percent of the full dataset is used to construct the Hessian. When there is no postfix, the Hessian matrix is formed exactly. 

\subsection{Binary Classification}
\label{sec:exp:binary}
\begin{figure}[htbp]
\begin{subfigure}{0.5\textwidth}
\includegraphics[width=1\linewidth]{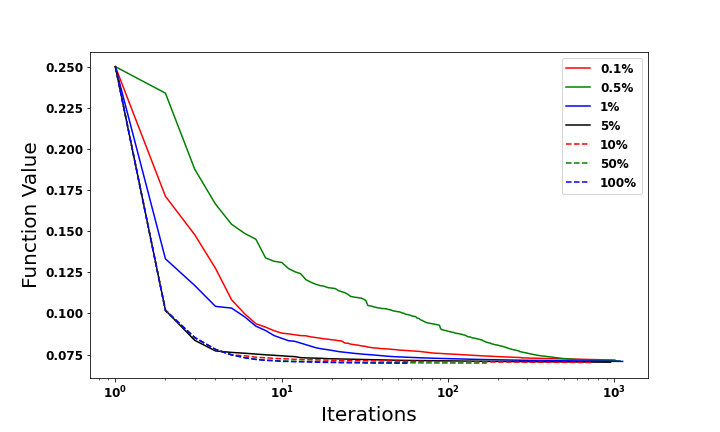} 
\end{subfigure}
\begin{subfigure}{0.5\textwidth}
\includegraphics[width=1\linewidth]{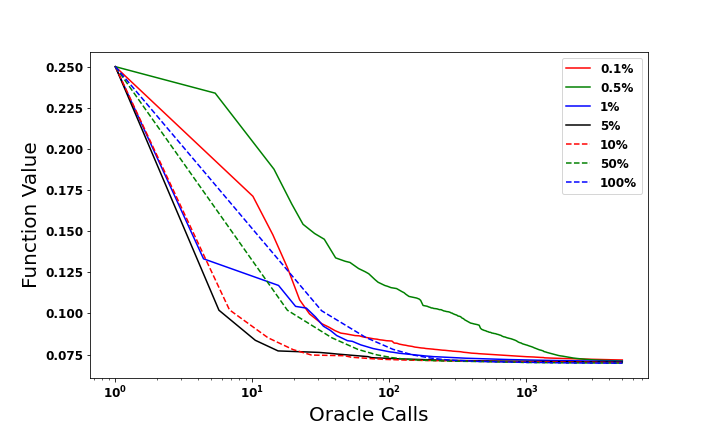}
\end{subfigure}
\caption{Performance of \cref{alg:NewtonMR} using various degrees of Hessian approximation on the nonconvex nonlinear least squares loss function. As predicted by our theory, \cref{alg:NewtonMR} converges irrespective of the degree of Hessian approximation. Also, Hessian approximation typically reduces computational costs; however, a substantial reduction in sub-sample size can lead to a significant loss of curvature information, resulting in poor performances. \label{fig:MNIST_logloss}}
\end{figure}

\begin{figure}[htbp]
\begin{subfigure}{0.5\textwidth}
\includegraphics[width=1\linewidth]{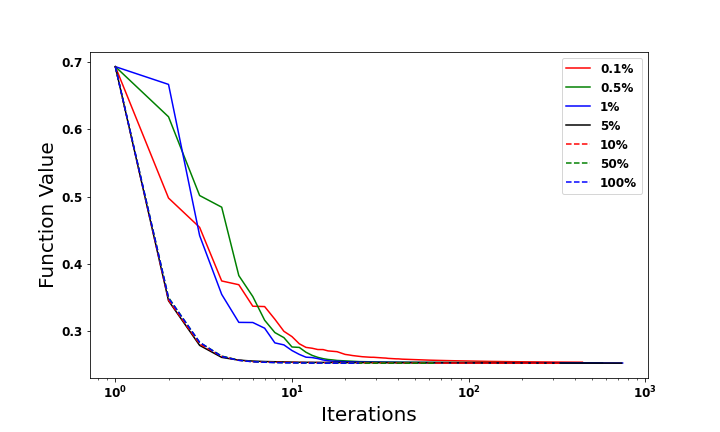} 
\end{subfigure}
\begin{subfigure}{0.5\textwidth}
\includegraphics[width=1\linewidth]{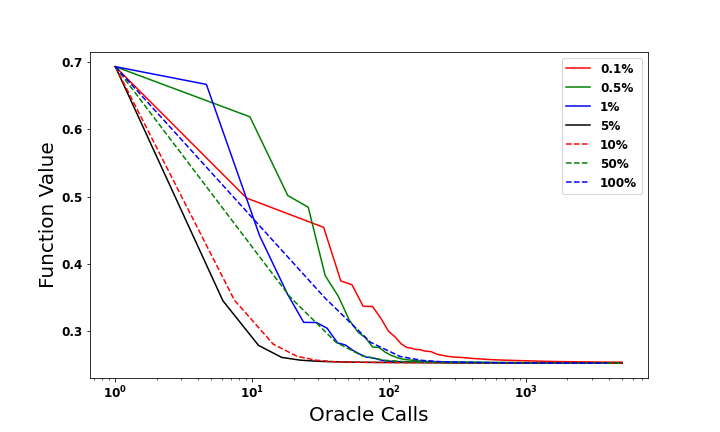}
\end{subfigure}
\caption{Performance of \cref{alg:NewtonMR} using various degrees of Hessian approximation on the convex logistic loss function.  As predicted by our theory, \cref{alg:NewtonMR} converges irrespective of the degree of Hessian approximation. Also, Hessian approximation typically reduces computational costs; however, a substantial reduction in sub-sample size can lead to a significant loss of curvature information, resulting in poor performances. \label{fig:MNIST_NLS}}
\end{figure}

We first investigate the advantages/disadvantages of Hessian approximation within Newton-MR framework using two simple finite-sum problems in the form of binary classification. Namely, we consider convex binary logistic regression (\cref{fig:MNIST_logloss}),
\begin{align*}
	f(\f{x}) = \frac{1}{n}\sum_{i=1}^n\ln(1 + \exp({\langle\f{a}_i,\f{x}\rangle})) - b_i\langle\f{a}_i,\f{x}\rangle,
\end{align*}
and nonconvex nonlinear least square (\cref{fig:MNIST_NLS}),
\begin{align*}
	f(\f{x}) = \frac{1}{n}\sum_{i=1}^n \left(\frac{1}{1 + \exp(-\langle \f{a}_i, \f{x} \rangle)} - b_i\right)^2.
\end{align*}
For our experiments here, we use the MINST data set \cite{deng2012mnist}, which contains of $60,000$ grayscale images of size $ 28 \times 28 $, i.e., $\{(\f{a}_i, b_i)\}_{i=1}^{n=60,000} \subset \mathbb{R}^{784}\times[0,1]$. The labels of MNIST dataset are converted to binary classification, i.e., even and odd numbers.

In \cref{alg:NewtonMR}, we set $\eta = 0.001$ and $\f{x}_{0} = \f{0}$. The results of running \cref{alg:NewtonMR} using various degrees of Hessian approximation are given in \cref{fig:MNIST_NLS,fig:MNIST_logloss}. As predicted by our theory, \cref{alg:NewtonMR} converges irrespective of the degree of Hessian approximation. Also, while Hessian approximation generally helps with reducing the overall computational costs,  its effectiveness can diminish if the approximation is overly crude. This is expected since a significant reduction in sub-sample size could lead to a substantial loss of curvature information, resulting in poor performance.

\subsection{Feed Forward Neural Network}
\label{sec:exp:ffnn}
In this section, we consider a multi-class classification problem using a feed forward neural network, namely 
\begin{equation}\label{eq:celoss}
f(\f{x}) = -\frac{1}{n}\left(\sum_{i=1}^n\dotprod{\f{b}_i,\ln\left(\f{h}(\f{x} ; \f{a}_i)\right)}\right) + \lambda g(\f{x}),
\end{equation}
where $\f{h}(\f{x} ; \f{a}_i) : \mathbb{R}^d \to [0,1]^{C}$ represent the output of the neural network, $ \f{b}_{i} \in \{0,1\}^{C} $ denotes the labels corresponding to input data $ \f{a}_{i} $, $g(\f{x}) : \mathbb{R}^d \to \mathbb{R}$ is a nonconvex regularizer defined as $g(\f{x}) \triangleq \sum^d_{i=1} x_i^2 / (1 + x_i^2)$, and $\lambda > 0$ is a regularization parameter. 

\begin{figure}[htbp]
\begin{subfigure}{0.5\textwidth}
\includegraphics[width=1\linewidth]{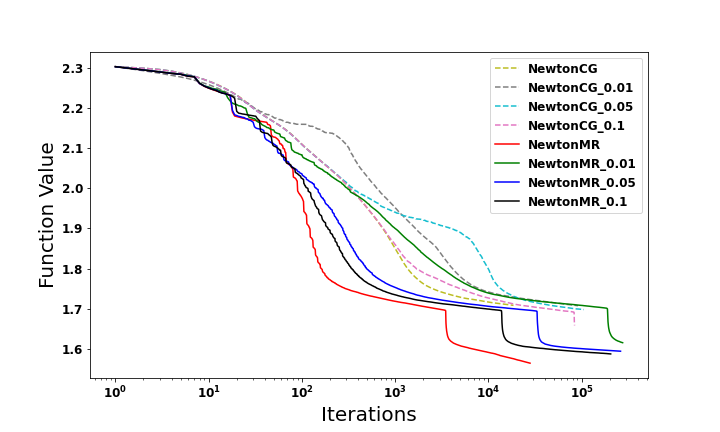} 
\end{subfigure}
\begin{subfigure}{0.5\textwidth}
\includegraphics[width=1\linewidth]{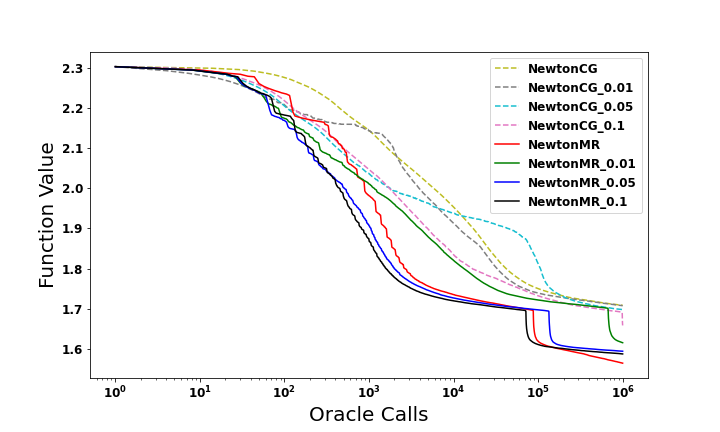}
\end{subfigure}
\caption{Comparison of Newton-MR and Newton-CG on CIFAR10 dataset in \cref{sec:exp:ffnn}. 
}\label{fig:cifar_cg}
\end{figure}

\begin{figure}[htbp]
\begin{subfigure}{0.5\textwidth}
\includegraphics[width=1\linewidth]{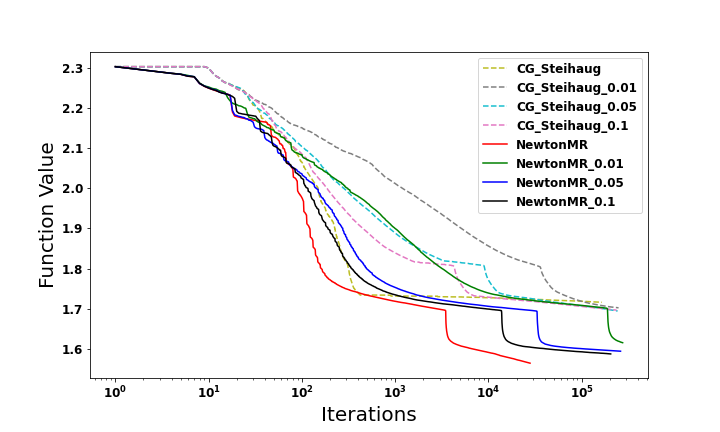} 
\end{subfigure}
\begin{subfigure}{0.5\textwidth}
\includegraphics[width=1\linewidth]{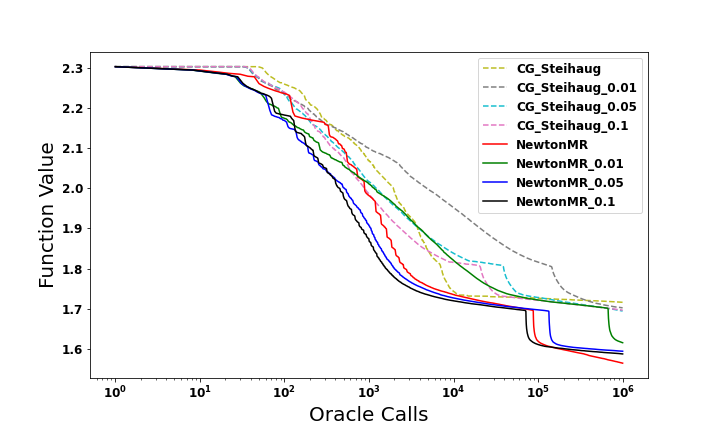} 
\end{subfigure}
\caption{Comparison of Newton-MR and Trust-Region on CIFAR10 dataset in \cref{sec:exp:ffnn}. 
}\label{fig:cifar_tr}
\end{figure}

\begin{figure}[htbp]
\begin{subfigure}{0.5\textwidth}
\includegraphics[width=1\linewidth]{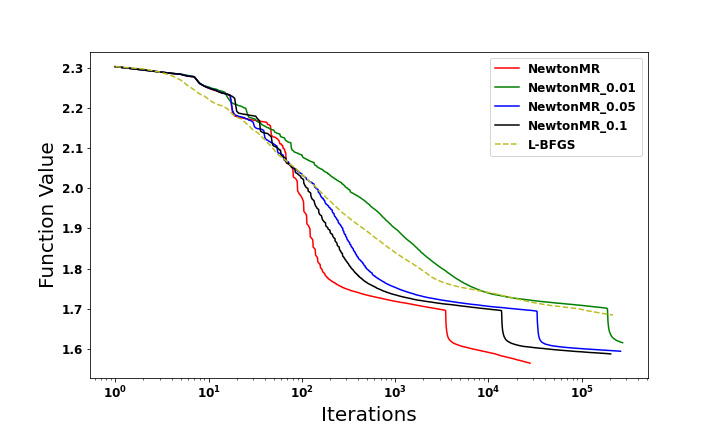} 
\end{subfigure}
\begin{subfigure}{0.5\textwidth}
\includegraphics[width=1\linewidth]{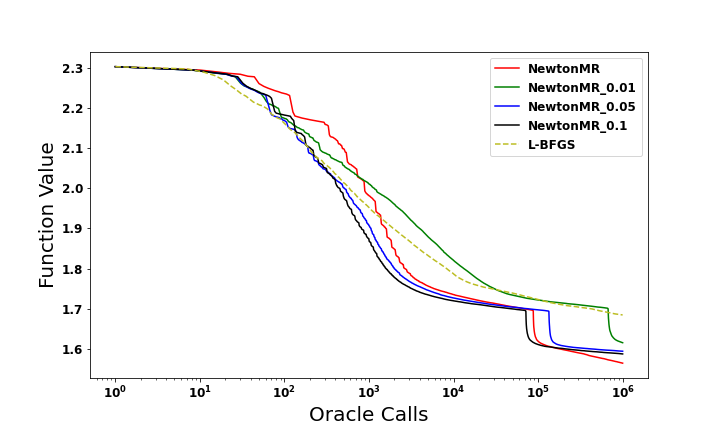}
\end{subfigure}
\caption{Comparison of Newton-MR and L-BFGS on CIFAR10 dataset in \cref{sec:exp:ffnn}. 
}\label{fig:cifar_lbfgs}
\end{figure}

We first consider CIFAR10 dataset \cite{krizhevsky2009learning}, consisting of $60,000$ colored images of size $ 32 \times 32 $ in 10 classes, i.e., $\{(\f{a}_i,\f{b}_i)\}_{i = 1}^{n = 60,000} \subset \mathbb{R}^{3,072} \times \{0, 1\}^{10}$. We split the 60,000 samples into 50,000 training samples and 10,000 validation samples. \cref{fig:cifar_cg,fig:cifar_tr,fig:cifar_lbfgs} depict the performance of each algorithm on the training data. Plots showing performance in terms of validation accuracy/error are gathered in \cref{appendix:validation}.
The feed forward neural network has the following architecture, 
\begin{align*}
	3072 \text{ (input)} \to \tanh \to 512 \to \tanh \to 128 \to \tanh \to 32 \to \text{softmax} \to 10 \text{ (output)},
\end{align*}
resulting in $d = 1,643,498$ numbers of parameters. The regularization parameter is set to $\lambda = 10^{-8}$. All algorithms are initialized from a normal distribution with zero mean and co-variance $0.1\f{I}$.
The inexactness and the $\sigma$-LC parameters of \cref{alg:NewtonMR} are set to $\eta = 10$ and $\sigma = 10^{-16}$. For Newton-CG \cite{royer2020newton,yao2023inexact}, the desired accuracy for Capped-CG is set to $0.01$. 

\begin{figure}[htbp]
\begin{subfigure}{0.5\textwidth}
\includegraphics[width=1\linewidth]{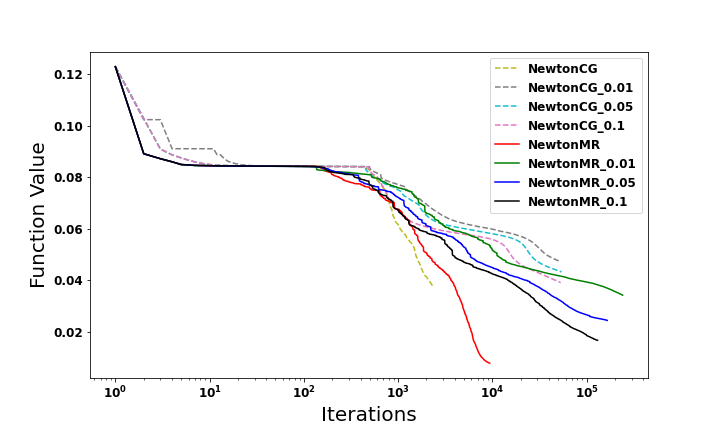} 
\end{subfigure}
\begin{subfigure}{0.5\textwidth}
\includegraphics[width=1\linewidth]{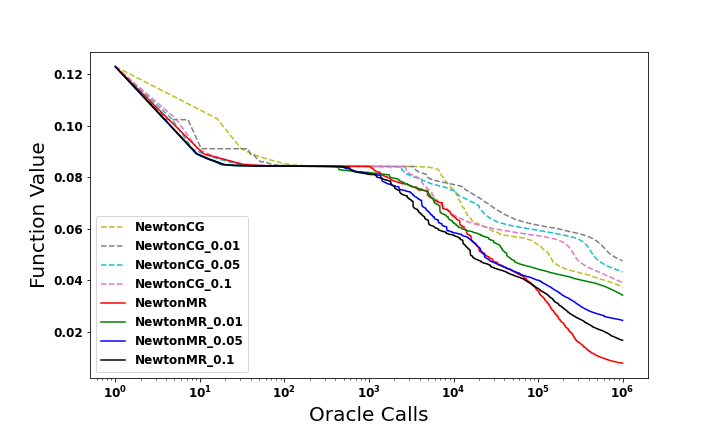}
\end{subfigure}
\caption{Comparison of Newton-MR and Newton-CG on Covertype dataset in \cref{sec:exp:ffnn}. 
\label{fig:cov_cg}}
\end{figure}

\begin{figure}[htbp]
\begin{subfigure}{0.5\textwidth}
\includegraphics[width=1\linewidth]{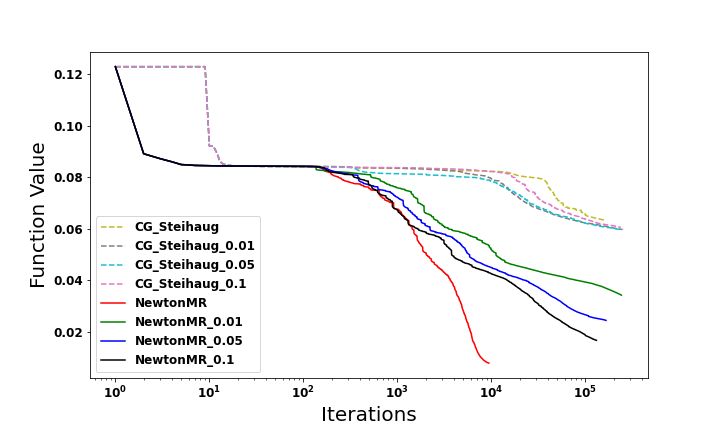} 
\end{subfigure}
\begin{subfigure}{0.5\textwidth}
\includegraphics[width=1\linewidth]{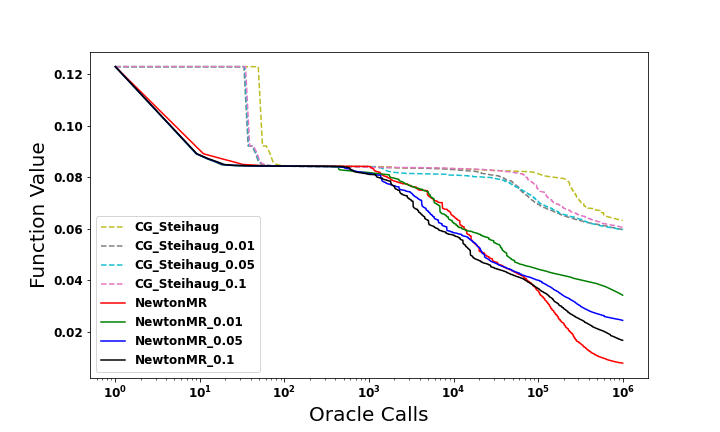}
\end{subfigure}
\caption{Comparison of Newton-MR and Trust-Region on Covertype dataset in \cref{sec:exp:ffnn}. 
\label{fig:cov_tr}}
\end{figure}

\begin{figure}[htbp]
\begin{subfigure}{0.5\textwidth}
\includegraphics[width=1\linewidth]{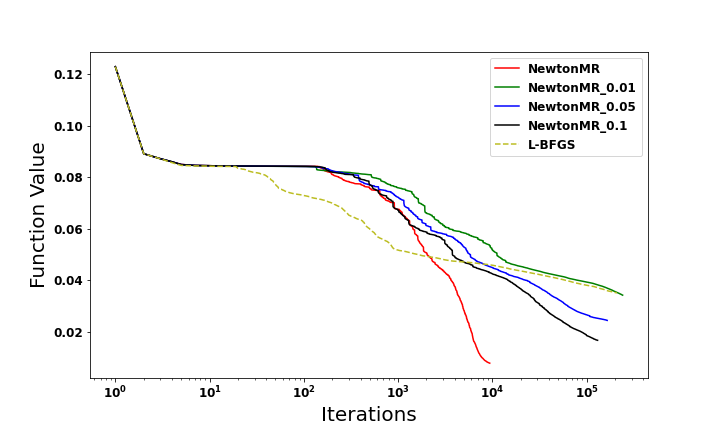} 
\end{subfigure}
\begin{subfigure}{0.5\textwidth}
\includegraphics[width=1\linewidth]{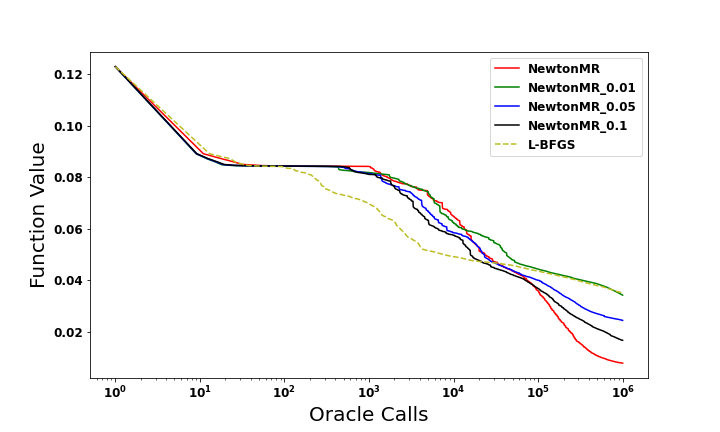}
\end{subfigure}
\caption{Comparison of Newton-MR and L-BFGS on Covertype dataset in \cref{sec:exp:ffnn}. 
\label{fig:cov_lbfgs}}
\end{figure}

Next, we consider a different dataset, namely Covertype \cite{misc_covertype_31}, which includes $n = 300,000$ training and $2,000$ validation samples of $54$ attributes across seven classes. \cref{fig:cov_tr,fig:cov_cg,fig:cov_lbfgs}  depict the performance of each algorithm on the training data. Plots showing performance in terms of validation accuracy/error are gathered in \cref{appendix:validation}. For these experiments, the neural network architecture is chosen to be 
\begin{align*}
	& 54 \text{ (input)} \to \sigmoid \to 256 \to \tanh \to 256 \to \sigmoid \to 128 \to \tanh \to 128 \to\\
    & \to \sigmoid \to 64 \to \tanh \to 64 \to \sigmoid \to 32 \to \tanh \to 32 \to \text{softmax} \to 7 \text{ (output)},
\end{align*}
resulting in $d = 145,063$. We also set $\lambda = 10^{-8}$. Again, $\f{x}_0$ is chosen as above. For Newton-MR, the sub-problem inexactness and $\sigma$-LC parameters are changed to $\eta = 10,000$ and $\sigma = 10^{-26}$ respectively. 
For Newton-CG, the desired accuracy parameter is set to $0.01$. 

In the experiments of this section, all algorithms reach the maximum allowed number of oracle calls. In all \cref{fig:cifar_cg,fig:cifar_tr,fig:cifar_lbfgs,fig:cov_tr,fig:cov_cg,fig:cov_lbfgs}, \cref{alg:NewtonMR} with appropriate degree of sub-sampling either is very competitive with, or downright superior to, alternative methods. Using the exact Hessian improves convergence by reducing the number of iterations. However, unless highly accurate solutions are required -- which is often not necessary in machine learning applications -- subsampling the Hessian can reduce the number of required oracle calls. What is striking is the relative poor performance of all Newton-CG variants. This further corroborates the observations in \cite{roosta2022newton,liu2021convergence,liu2022newton} that MINRES typically manifests itself as a far superior inner solver than the widely used CG method. 

\subsection{Recurrent Neural Network}
\label{sec:exp:rnn}
Finally, we consider a slightly more challenging problem of training a recurrent neural network (RNN) \cite{sherstinsky2020fundamentals} with mean squared error loss. Specifically, we consider the following objective function 
\begin{align*}
	f(\f{x}) = \frac{1}{n}\sum_{i=1}^n \left\|\f{h}(\f{x};\f{a}_i) - \f{b}_i\right\|^{2} + \lambda g(\f{x}),
\end{align*}
where $ \f{h}(\f{x};\f{a}) $ is the output of the RNN and $g(\f{x})$ is again the same nonconvex regularizer as in \cref{sec:exp:ffnn} with regularization parameter $\lambda = 10^{-8}$. The point $\f{x}_0$ is chosen in the same manner as in \cref{sec:exp:ffnn}.

For the experiment here, we use the Gas sensor array under dynamic gas mixtures data set \cite{fonollosa2015reservoir}, which is a time series data. In particular, the data contains a time series of 16 features across three time stamps. So, we get $\{\f{a}_i, \f{b}_i\}_{i = 1}^{n=45,828} \subset (\mathbb{R}^3 \times \mathbb{R}^{16}) \times \mathbb{R}$, and $ d = 100,417$. We further split the samples into 40,828 training and 5,000 validation samples. For Newton-MR, the relative residual tolerance is set to $\eta = 100$. For Newton-CG, the accuracy of Capped-CG is set to $0.1$.

\begin{figure}[htbp]
\begin{subfigure}{0.5\textwidth}
\includegraphics[width=1\linewidth]{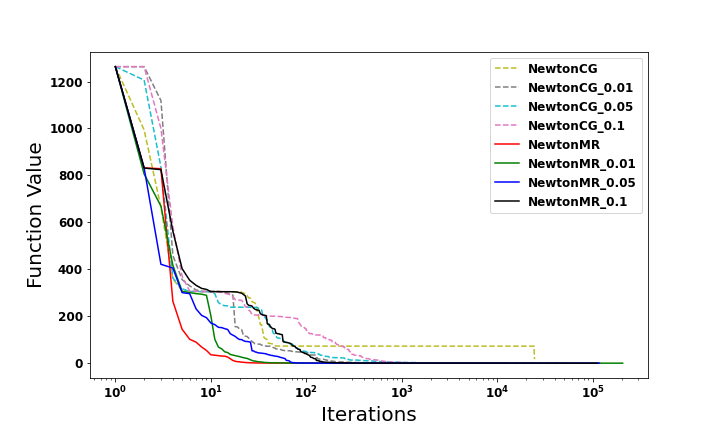} 
\end{subfigure}
\begin{subfigure}{0.5\textwidth}
\includegraphics[width=1\linewidth]{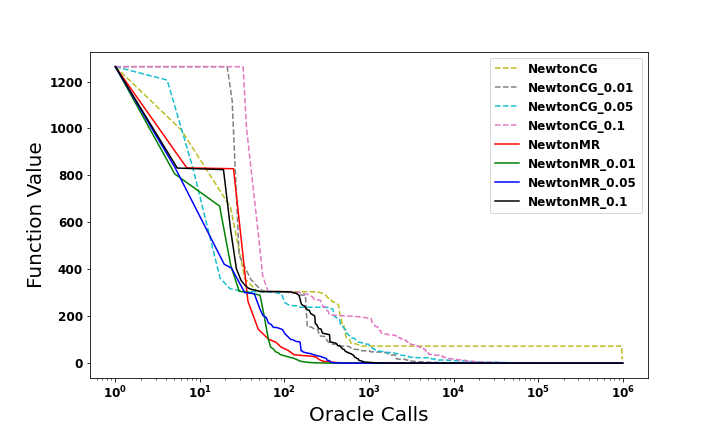}
\end{subfigure}
\caption{Comparison of Newton-MR and Newton-CG for training a recurrent neural network as in \cref{sec:exp:rnn}. 
$\texttt{NewtonCG}\_0.05$ and $\texttt{NewtonCG}\_0.1$ are terminated after 1410 and 2168 iterations with function values of 2.4428 and 0.31219, respectively, as their step sizes get below $10^{-18}$.
\label{fig:eth_cg}}
\end{figure}

\begin{figure}[htbp]
\begin{subfigure}{0.5\textwidth}
\includegraphics[width=1\linewidth]{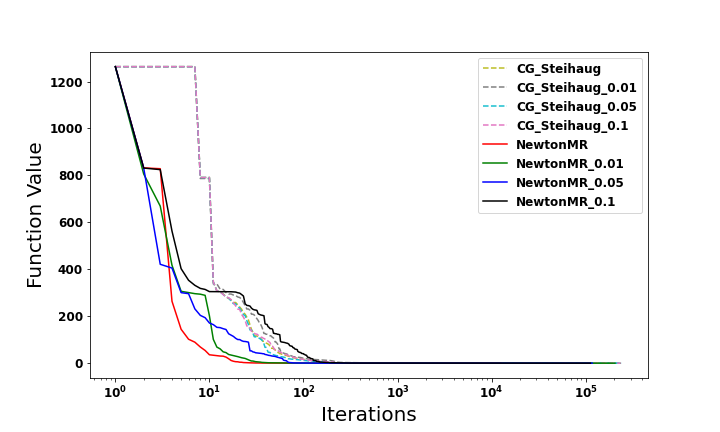} 
\end{subfigure}
\begin{subfigure}{0.5\textwidth}
\includegraphics[width=1\linewidth]{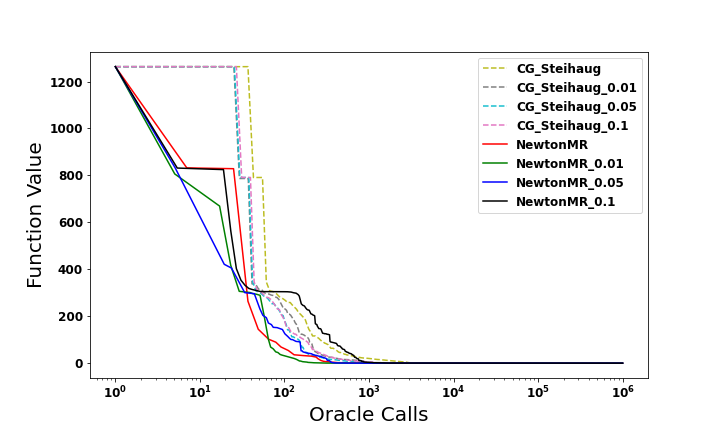}
\end{subfigure}
\caption{Comparison of Newton-MR and Trust-Region for training a recurrent neural network as in \cref{sec:exp:rnn}. 
\label{fig:eth_tr}}
\end{figure}

\begin{figure}[htbp]
\begin{subfigure}{0.5\textwidth}
\includegraphics[width=1\linewidth]{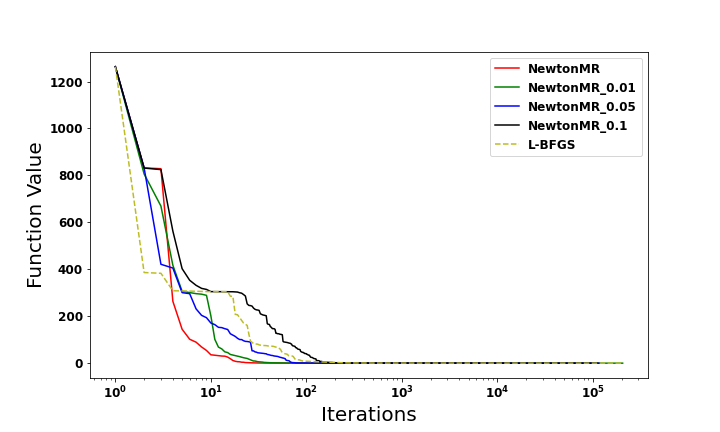} 
\end{subfigure}
\begin{subfigure}{0.5\textwidth}
\includegraphics[width=1\linewidth]{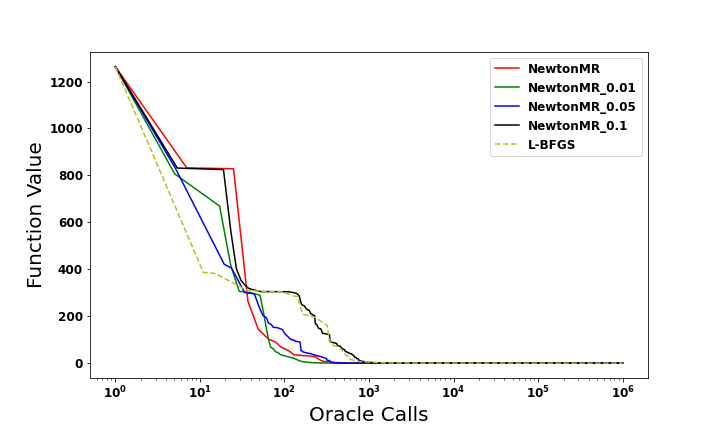}
\end{subfigure}
\caption{Comparison of Newton-MR and L-BFGS  for training a recurrent neural network as in \cref{sec:exp:rnn}. 
\label{fig:eth_lbfgs}}
\end{figure}

\cref{fig:eth_cg,fig:eth_lbfgs,fig:eth_tr} depict the performance of each algorithm on the training data. Plots showing performance in terms of validation accuracy/error are gathered in \cref{appendix:validation}.
Once again, in all these experiments, \cref{alg:NewtonMR} with adequate sub-sampling can either outperform or at least remain competitive with the alternative methods. Also, one can again observe the visibly poor performance of Newton-CG relative to \cref{alg:NewtonMR}. 
In \cref{fig:eth_cg}, $\texttt{NewtonCG}\_0.05$ and $\texttt{NewtonCG}\_0.1$ are terminated after 1410 and 2168 iterations, respectively, as their step sizes get below $10^{-18}$. All the other algorithms reach the maximum allowed number of oracle calls.

\section{Conclusion}\label{sec:conclusion}
We considered a variant of the Newton-MR algorithm for nonconvex problems, initially proposed in \cite{liu2021convergence}, to accommodate inexact Hessian information. 
Our primary focus was to establish iteration and operation complexities of our proposed variant under simple assumptions and with minimal algorithmic changes. In particular, we assume minimal smoothness on the gradient without extending it to the Hessian, and we avoid Hessian regularization/damping techniques to maintain curvature information and prevent performance degradation. We show that our algorithm converges regardless of the degree of approximation of the Hessian as well as the accuracy of the solution to the sub-problem. Beyond general nonconvex settings, we consider, in particular, a subclass of functions that satisfy the $\theta$-PL condition and show that with $\theta = 2$, our algorithm achieves a global linear convergence rate while for $1 \leq \theta < 2$, the global convergence transitions from a linear rate to a sub-linear rate in the neighbourhood of a solution. Finally, we provided numerical experiments on several machine learning problems to further illustrate the advantages of sub-sampled Hessian in the Newton-MR framework.

\subsection{Acknowledgments}
Fred Roosta was partially supported by the Australian Research Council through an Industrial Transformation Training Centre for Information Resilience (IC200100022) as well as a Discovery Early Career Researcher Award (DE180100923).

\bibliographystyle{plain}
\bibliography{reference}

\appendix
\section{Appendix}\label{sec:appendix}
\subsection{Proof of \cref{property:MINRES}}

The results have been established, in one form or another, in prior works \cite{saad2003iterative,liu2022minres}; however we provide the proofs here for the sake of completeness. The proofs of \cref{prop:monotonicity,prop:residual_descent} can be readily found in \cite{saad2003iterative} and \cite[Lemma 3.1d]{liu2022minres} respectively. 
For \cref{prop:inner_residual}, we note $\langle\f{r}^{(t-1)},\f{r}^{(i)}\rangle = \langle\f{r}^{(t-1)},-\f{g} - \inH\f{s}^{(i)}\rangle = \|\f{r}^{(t-1)}\|^2$, where the last equality follows from \cref{prop:residual_descent} and Petrov-Galerkin condition, i.e., $\f{r}^{(t-1)} \perp \inH\krylov{t}{\inH}{\f{g}}$. 
Property \cref{prop:descent_in_quadratic} can be found in \cite[Theorem 3.8a]{liu2022minres}. 
For \cref{prop:increasing_s}, the increasing sequence $\|\f{s}^{(i)}\| \geq \|\f{s}^{(j)}\|$ can be again found in \cite[Theorem 3.11e]{liu2022minres}. 
For the expression of $\f{s}^{(1)}$, we write
\begin{align*}
    \f{s}^{(1)} = \underset{\f{s}\in\krylov{1}{\inH}{\f{g}}}{\text{argmin}}\vnorm{-\f{g} - \inH\f{s}}^2 = \underset{\f{s}\in\text{span}(\f{g})}{\text{argmin}}\vnorm{\inH\f{s}}^2 + 2\lr{\f{g},\inH\f{s}}
\end{align*}
So, the vector $\f{s}^{(1)}$ must be some multiple of $\f{g}$, i.e., $\f{s}^{(1)} = \alpha\f{g}$. To find $\alpha$, we note,
\begin{align*}
    \frac{\lr{\f{g},\inH\f{g}}}{\vnorm{\inH\f{g}}^2} = \underset{\alpha\in\mathbb{R}}{\text{argmin}} \,\, \alpha^2\vnorm{\inH\f{g}}^2 + 2\alpha\lr{\f{g},\inH\f{g}}.
\end{align*}

\subsection{Further Numerical Results: Training and Validation Accuracy/Error}
\label{appendix:validation}
In this section, we provide the comparative performance of all methods, as measured by training and validation accuracy/error in each experiment of \cref{sec:numerical_exp}.

\begin{figure}[htbp]
\begin{subfigure}{0.5\textwidth}
\includegraphics[width=1\linewidth]{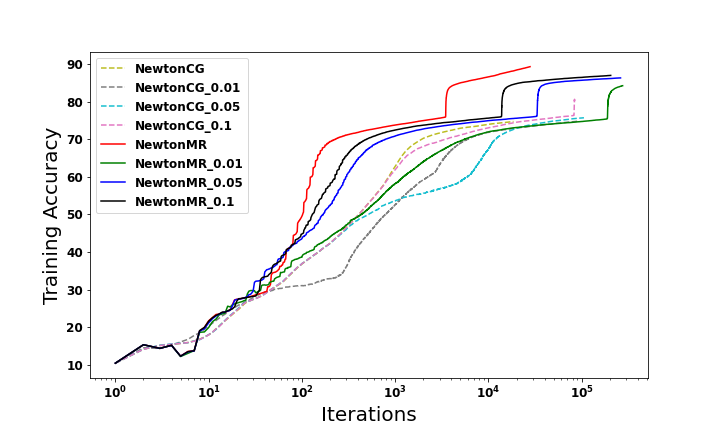} 
\end{subfigure}
\begin{subfigure}{0.5\textwidth}
\includegraphics[width=1\linewidth]{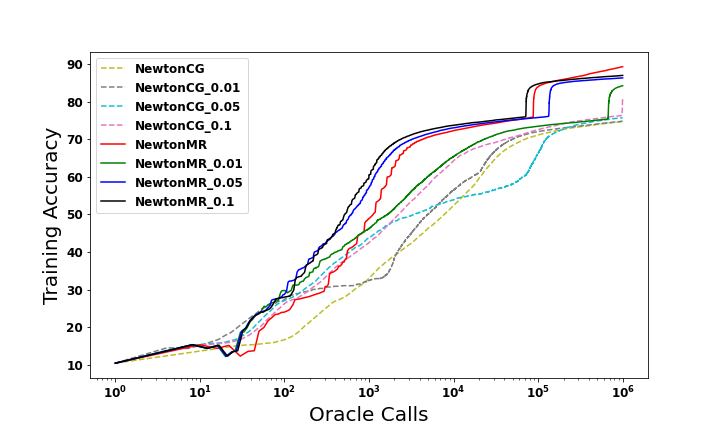}
\end{subfigure}
\begin{subfigure}{0.5\textwidth}
\includegraphics[width=1\linewidth]{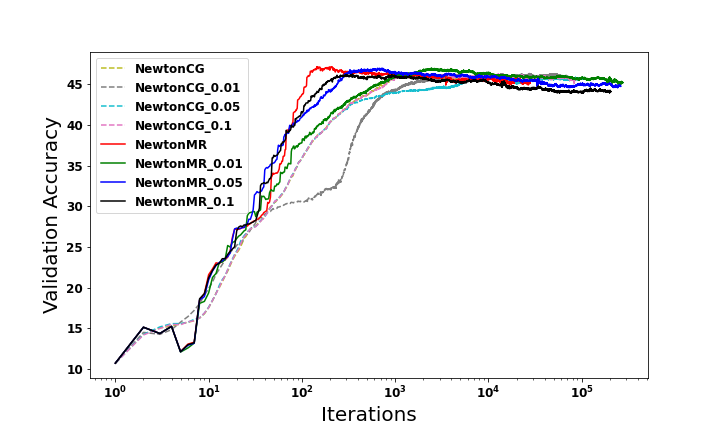} 
\end{subfigure}
\begin{subfigure}{0.5\textwidth}
\includegraphics[width=1\linewidth]{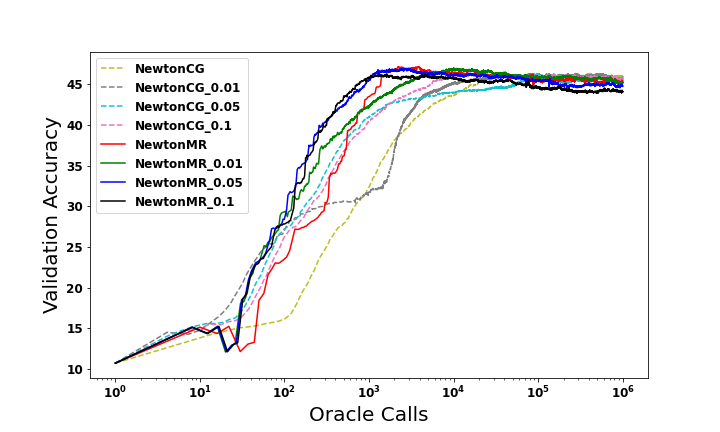}
\end{subfigure}
\caption{Comparison of Newton-MR and Newton-CG on CIFAR10 dataset in \cref{sec:exp:ffnn}.}\label{fig:further_cifar_cg}
\end{figure}

\begin{figure}[htbp]
\begin{subfigure}{0.5\textwidth}
\includegraphics[width=1\linewidth]{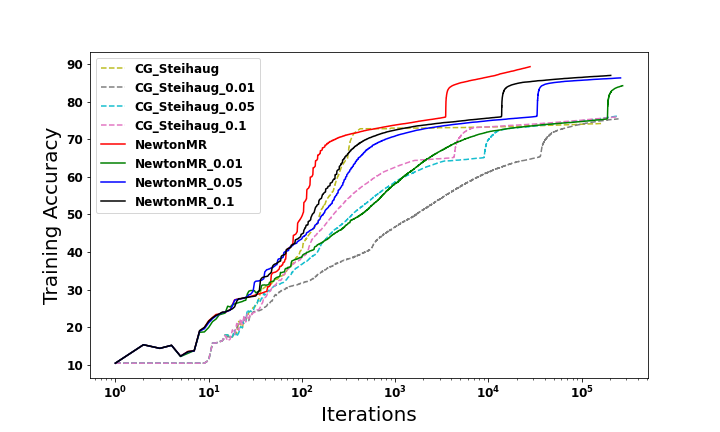} 
\end{subfigure}
\begin{subfigure}{0.5\textwidth}
\includegraphics[width=1\linewidth]{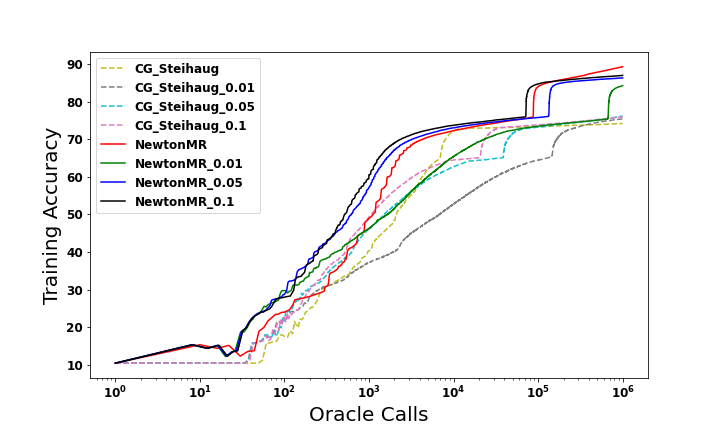}
\end{subfigure}
\begin{subfigure}{0.5\textwidth}
\includegraphics[width=1\linewidth]{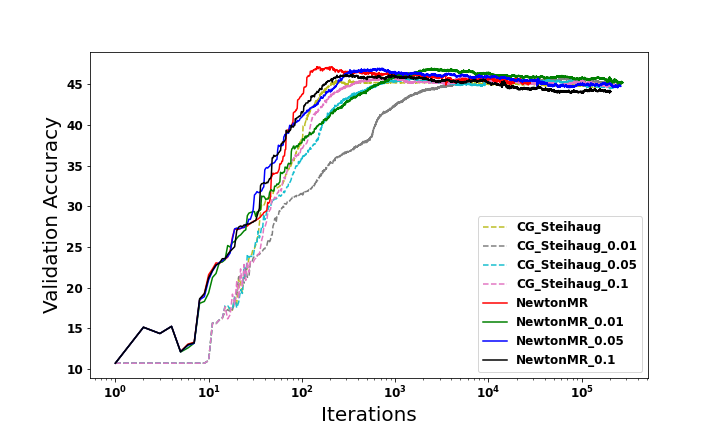} 
\end{subfigure}
\begin{subfigure}{0.5\textwidth}
\includegraphics[width=1\linewidth]{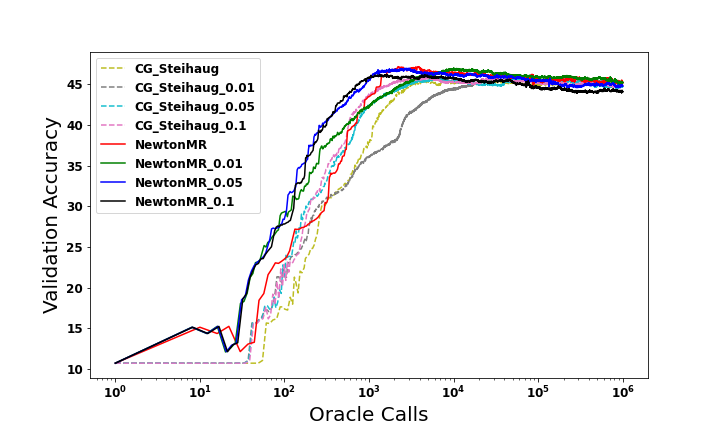}
\end{subfigure}
\caption{Comparison of Newton-MR and Trust-Region on CIFAR10 dataset in \cref{sec:exp:ffnn}.}\label{fig:further_cifar_tr}
\end{figure}

\begin{figure}[htbp]
\begin{subfigure}{0.5\textwidth}
\includegraphics[width=1\linewidth]{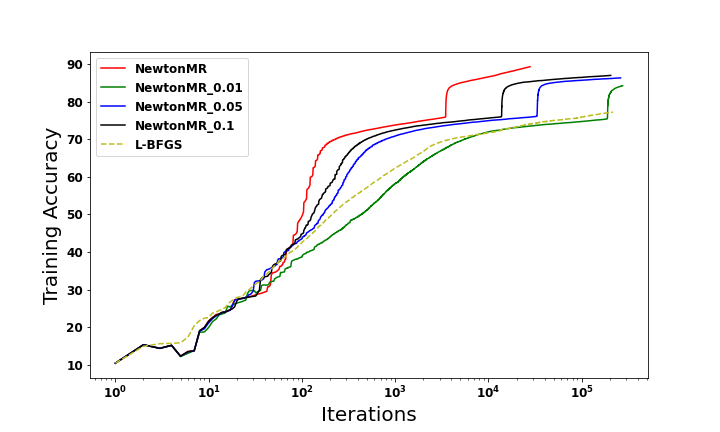} 
\end{subfigure}
\begin{subfigure}{0.5\textwidth}
\includegraphics[width=1\linewidth]{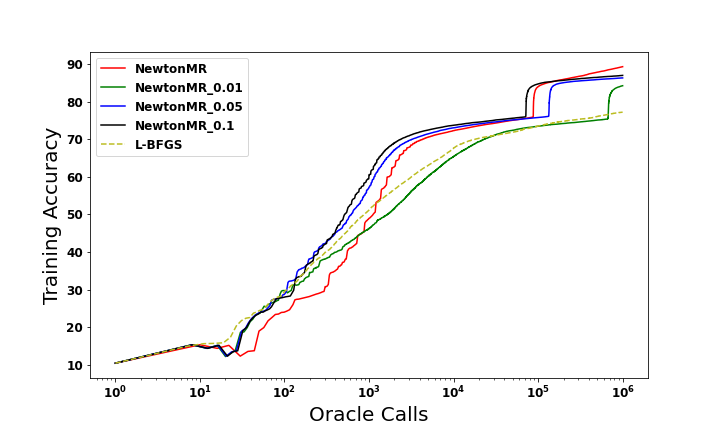}
\end{subfigure}
\begin{subfigure}{0.5\textwidth}
\includegraphics[width=1\linewidth]{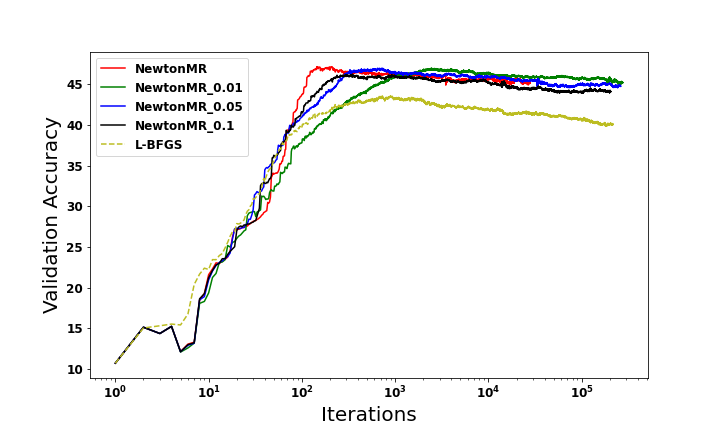} 
\end{subfigure}
\begin{subfigure}{0.5\textwidth}
\includegraphics[width=1\linewidth]{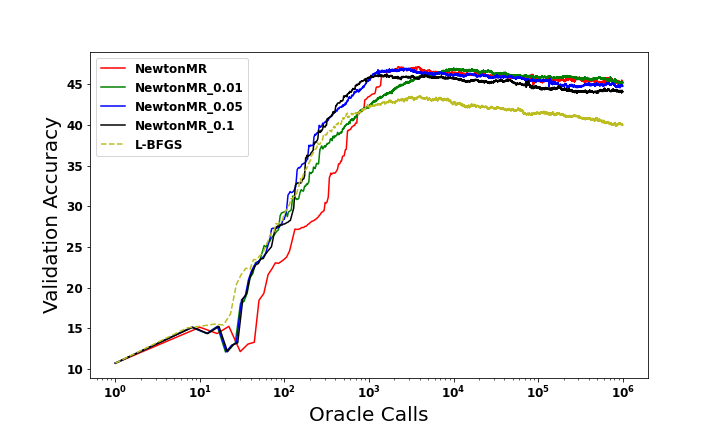}
\end{subfigure}
\caption{Comparison of Newton-MR and L-BFGS on CIFAR10 dataset in \cref{sec:exp:ffnn}.}\label{fig:further_cifar_lbfgs}
\end{figure}

\begin{figure}[htbp]
\begin{subfigure}{0.5\textwidth}
\includegraphics[width=1\linewidth]{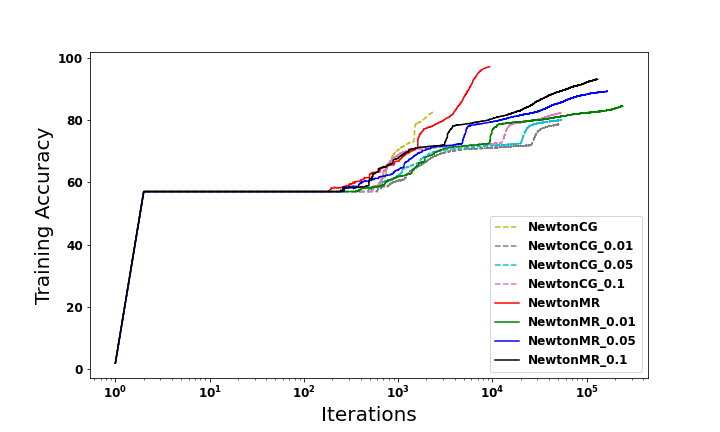} 
\end{subfigure}
\begin{subfigure}{0.5\textwidth}
\includegraphics[width=1\linewidth]{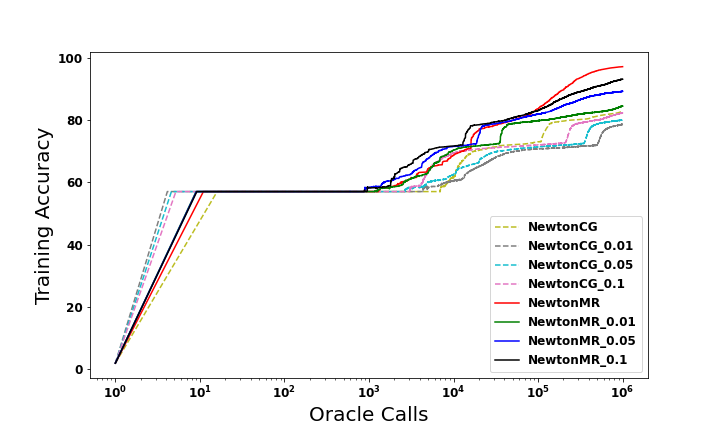}
\end{subfigure}
\begin{subfigure}{0.5\textwidth}
\includegraphics[width=1\linewidth]{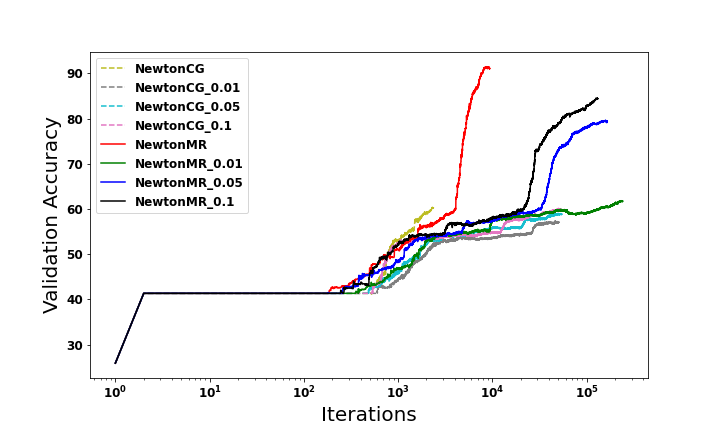} 
\end{subfigure}
\begin{subfigure}{0.5\textwidth}
\includegraphics[width=1\linewidth]{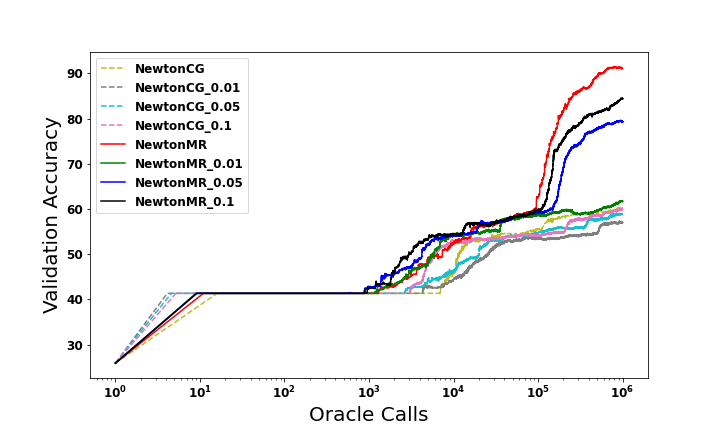}
\end{subfigure}
\caption{Comparison of Newton-MR and Newton-CG on Covertype dataset in \cref{sec:exp:ffnn}. 
\label{fig:further_cov_cg}}
\end{figure}

\begin{figure}[htbp]
\begin{subfigure}{0.5\textwidth}
\includegraphics[width=1\linewidth]{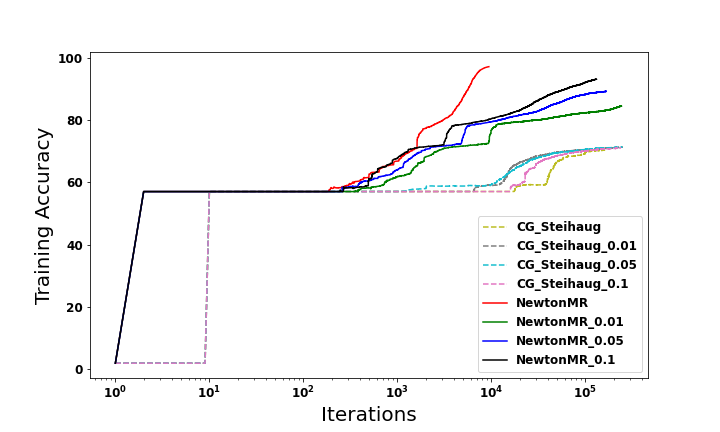} 
\end{subfigure}
\begin{subfigure}{0.5\textwidth}
\includegraphics[width=1\linewidth]{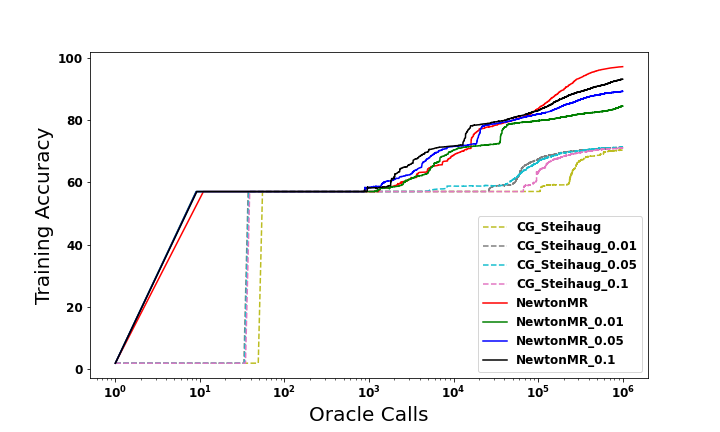}
\end{subfigure}
\begin{subfigure}{0.5\textwidth}
\includegraphics[width=1\linewidth]{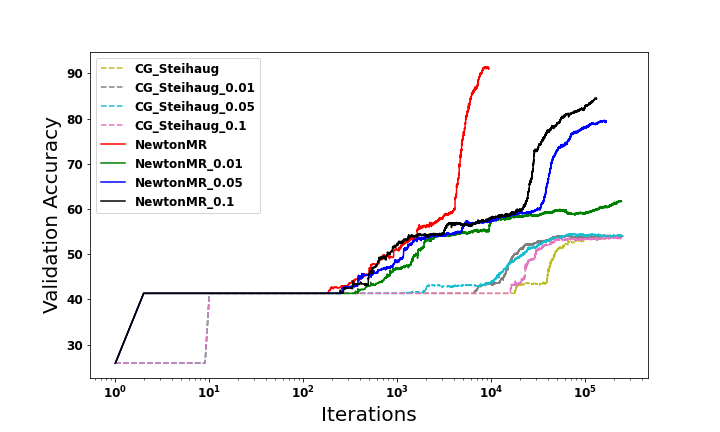} 
\end{subfigure}
\begin{subfigure}{0.5\textwidth}
\includegraphics[width=1\linewidth]{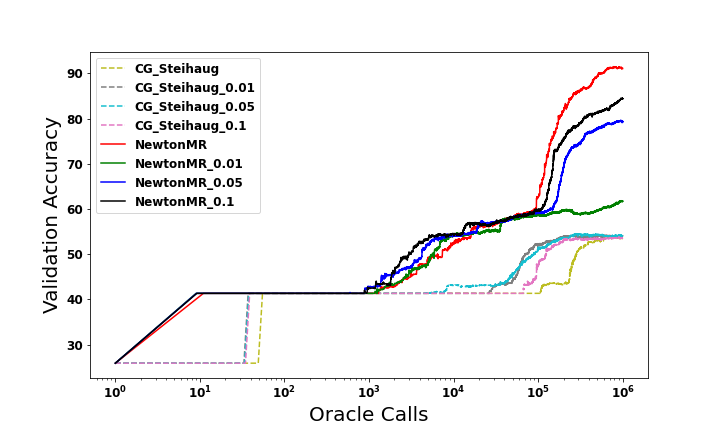}
\end{subfigure}
\caption{Comparison of Newton-MR and Trust-Region on Covertype dataset in \cref{sec:exp:ffnn}. 
\label{fig:further_cov_tr}}
\end{figure}

\begin{figure}[htbp]
\begin{subfigure}{0.5\textwidth}
\includegraphics[width=1\linewidth]{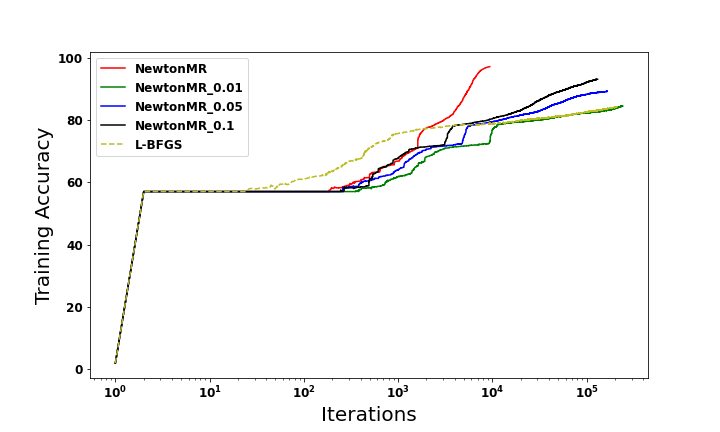} 
\end{subfigure}
\begin{subfigure}{0.5\textwidth}
\includegraphics[width=1\linewidth]{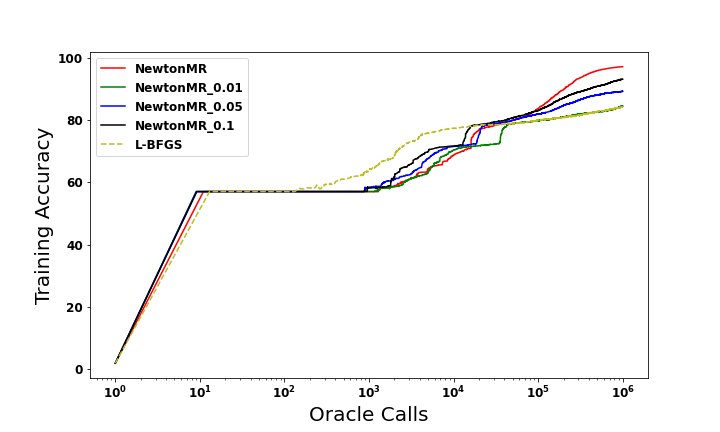}
\end{subfigure}
\begin{subfigure}{0.5\textwidth}
\includegraphics[width=1\linewidth]{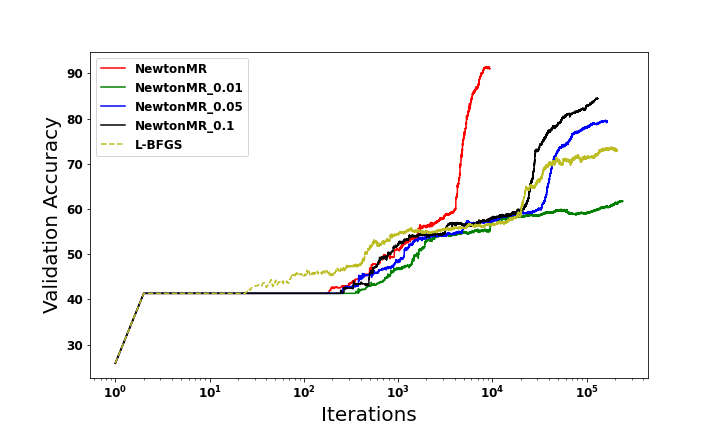} 
\end{subfigure}
\begin{subfigure}{0.5\textwidth}
\includegraphics[width=1\linewidth]{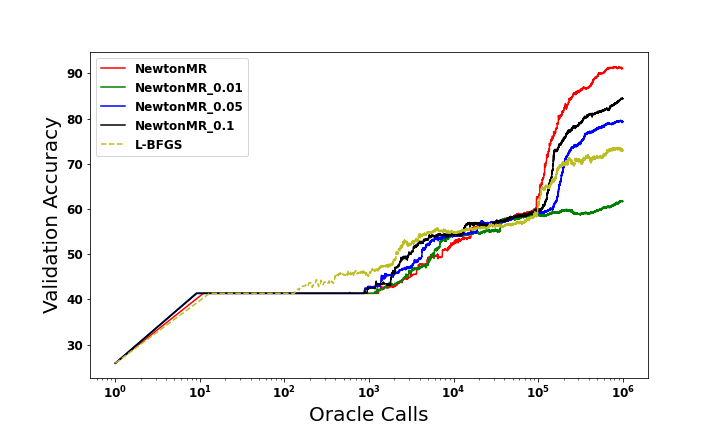}
\end{subfigure}
\caption{Comparison of Newton-MR and L-BFGS on Covertype in \cref{sec:exp:ffnn}. 
\label{fig:further_cov_lbfgs}}
\end{figure}

\begin{figure}[htbp]
\begin{subfigure}{0.5\textwidth}
\includegraphics[width=1\linewidth]{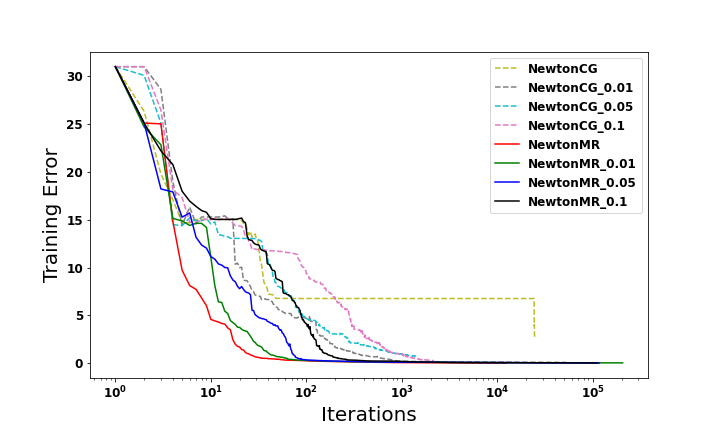} 
\end{subfigure}
\begin{subfigure}{0.5\textwidth}
\includegraphics[width=1\linewidth]{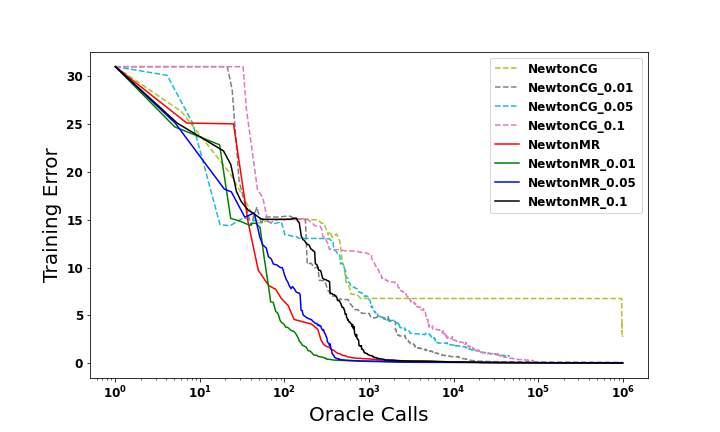}
\end{subfigure}
\begin{subfigure}{0.5\textwidth}
\includegraphics[width=1\linewidth]{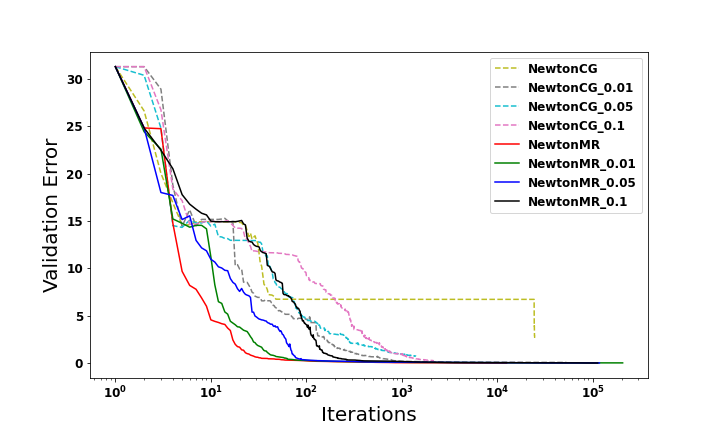} 
\end{subfigure}
\begin{subfigure}{0.5\textwidth}
\includegraphics[width=1\linewidth]{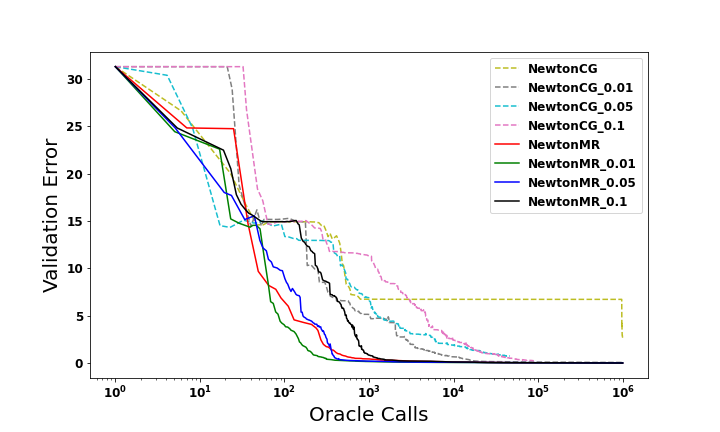}
\end{subfigure}
\caption{Comparison of Newton-MR and Newton-CG for training a recurrent neural network as in \cref{sec:exp:rnn}. The error is measured as $(\sum_{i \in \mathcal{S}} |h(\f{x}; \mathbf{a}_i) - b_i|) / \mathcal{|S|}$, where $\mathcal{S}$ denotes the training or validation set, containing $|\mathcal{S}|$ data points. Note that although the plots of validation error and training error appear highly similar, there are slight, albeit imperceptible, differences between them.
\label{fig:further_rnn_cg}}
\end{figure}

\begin{figure}[htbp]
\begin{subfigure}{0.5\textwidth}
\includegraphics[width=1\linewidth]{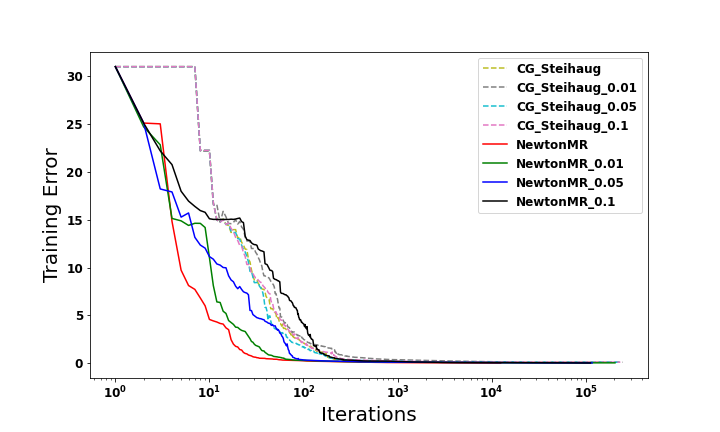} 
\end{subfigure}
\begin{subfigure}{0.5\textwidth}
\includegraphics[width=1\linewidth]{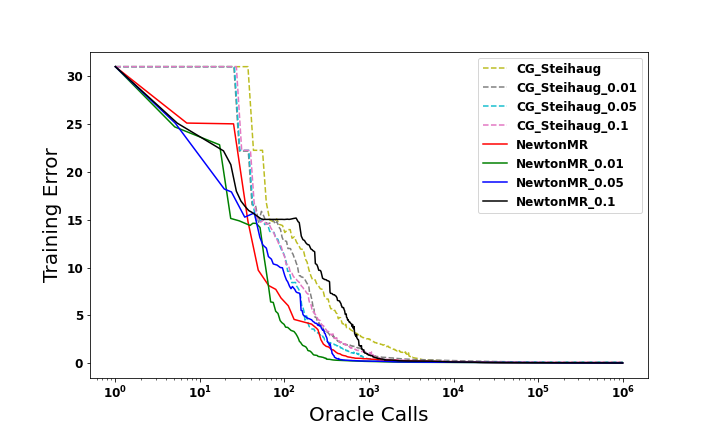}
\end{subfigure}
\begin{subfigure}{0.5\textwidth}
\includegraphics[width=1\linewidth]{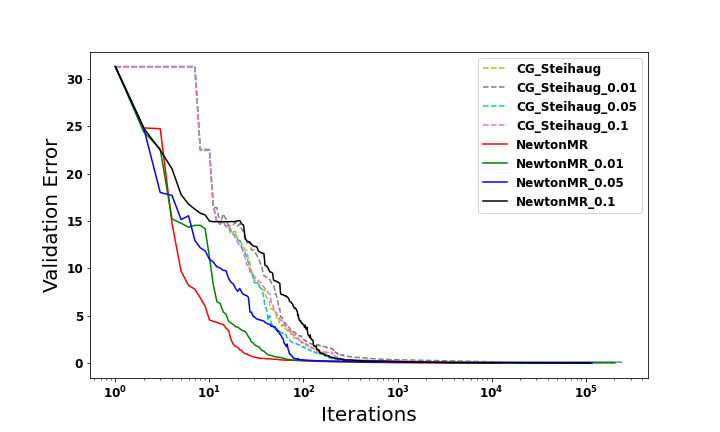} 
\end{subfigure}
\begin{subfigure}{0.5\textwidth}
\includegraphics[width=1\linewidth]{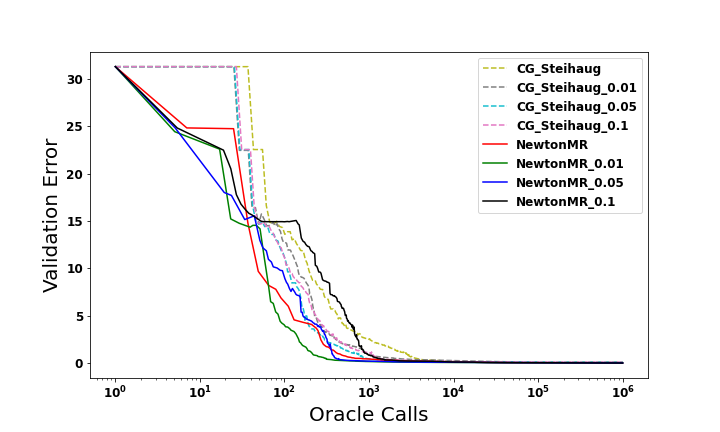}
\end{subfigure}
\caption{Comparison of Newton-MR and Trust-Region for training a recurrent neural network as in \cref{sec:exp:rnn}. The error is measured as $(\sum_{i \in \mathcal{S}} |h(\f{x}; \mathbf{a}_i) - b_i|) / \mathcal{|S|}$, where $\mathcal{S}$ denotes the training or validation set, containing $|\mathcal{S}|$ data points. Note that although the plots of validation error and training error appear highly similar, there are slight, albeit imperceptible, differences between them.
\label{fig:further_rnn_tr}}
\end{figure}

\begin{figure}[htbp]
\begin{subfigure}{0.5\textwidth}
\includegraphics[width=1\linewidth]{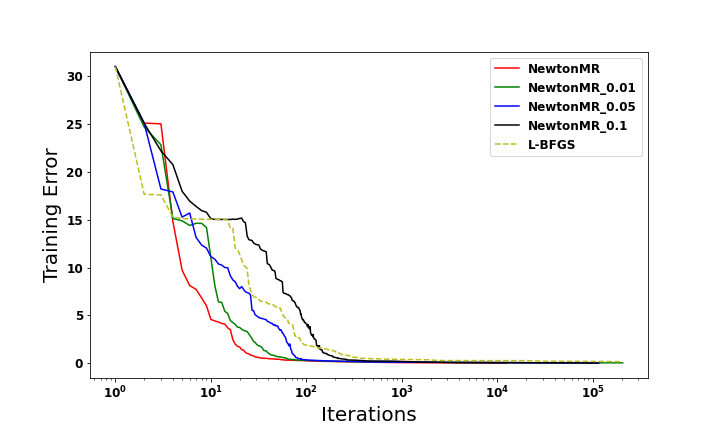} 
\end{subfigure}
\begin{subfigure}{0.5\textwidth}
\includegraphics[width=1\linewidth]{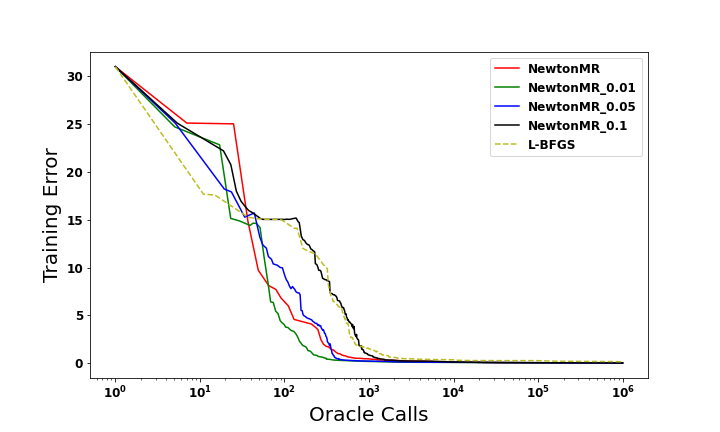}
\end{subfigure}
\begin{subfigure}{0.5\textwidth}
\includegraphics[width=1\linewidth]{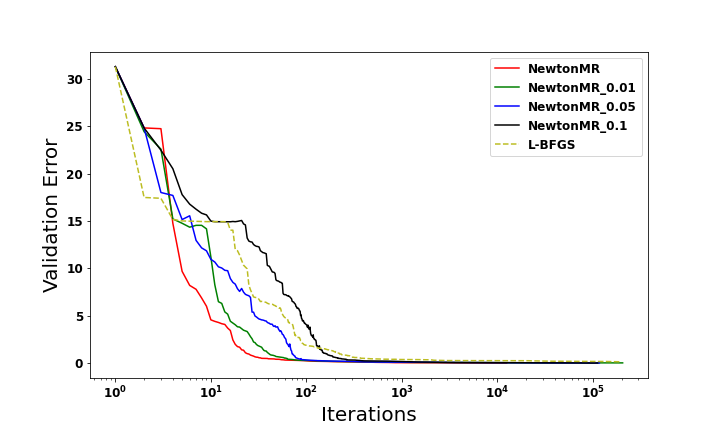} 
\end{subfigure}
\begin{subfigure}{0.5\textwidth}
\includegraphics[width=1\linewidth]{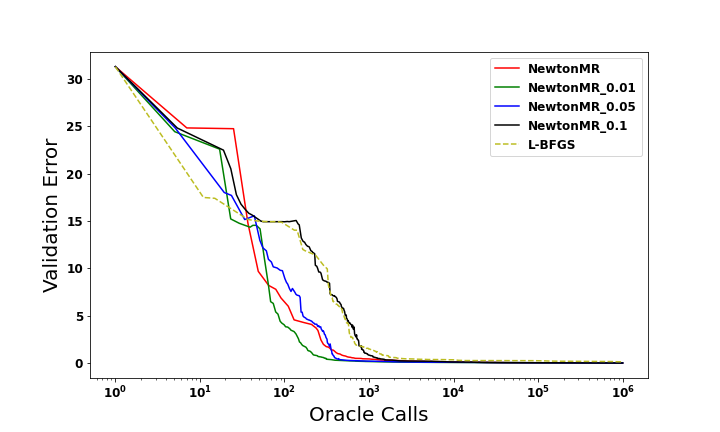}
\end{subfigure}
\caption{Comparison of Newton-MR and L-BFGS for training a recurrent neural network as in \cref{sec:exp:rnn}. The error is measured as $(\sum_{i \in \mathcal{S}} |h(\f{x}; \mathbf{a}_i) - b_i|) / \mathcal{|S|}$, where $\mathcal{S}$ denotes the training or validation set, containing $|\mathcal{S}|$ data points. Note that although the plots of validation error and training error appear highly similar, there are slight, albeit imperceptible, differences between them.
\label{fig:further_rnn_lbfgs}}
\end{figure}

\subsection{Why not damp/regularize the Hessian matrix?}
\label{sec:hessiam_damp}
\begin{figure}[htbp]
\begin{subfigure}{0.33\textwidth}
\includegraphics[width=1\linewidth]{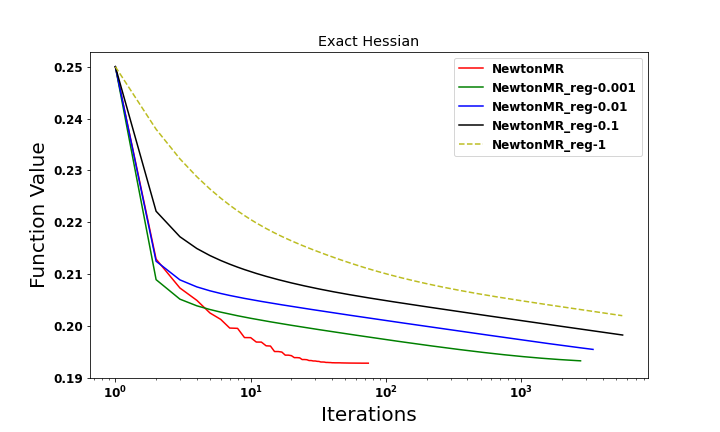} 
\end{subfigure}
\begin{subfigure}{0.33\textwidth}
\includegraphics[width=1\linewidth]{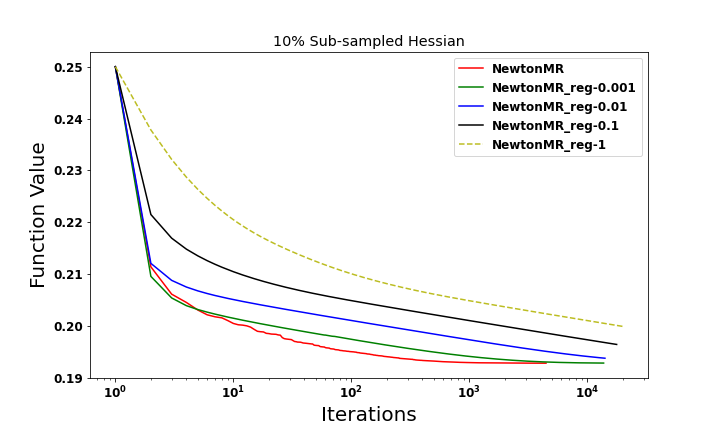} 
\end{subfigure}
\begin{subfigure}{0.33\textwidth}
\includegraphics[width=1\linewidth]{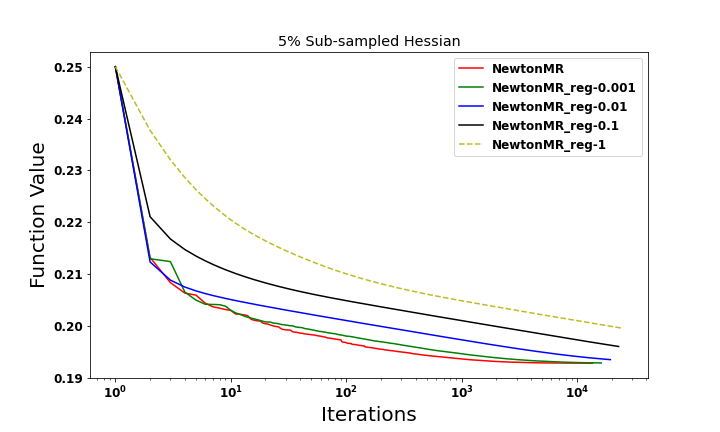} 
\end{subfigure}
\begin{subfigure}{0.33\textwidth}
\includegraphics[width=1\linewidth]{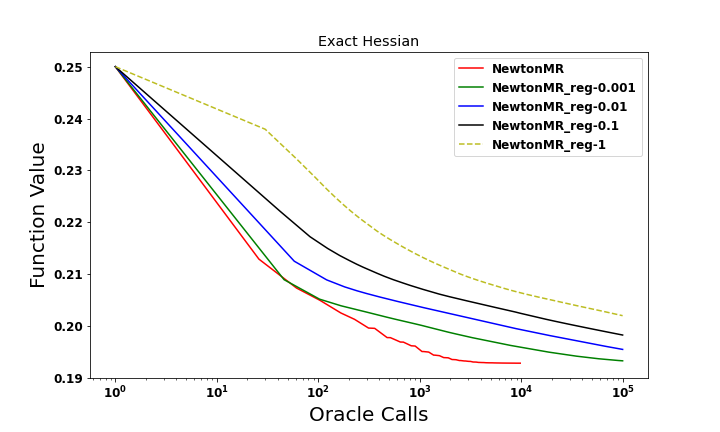} 
\end{subfigure}
\begin{subfigure}{0.33\textwidth}
\includegraphics[width=1\linewidth]{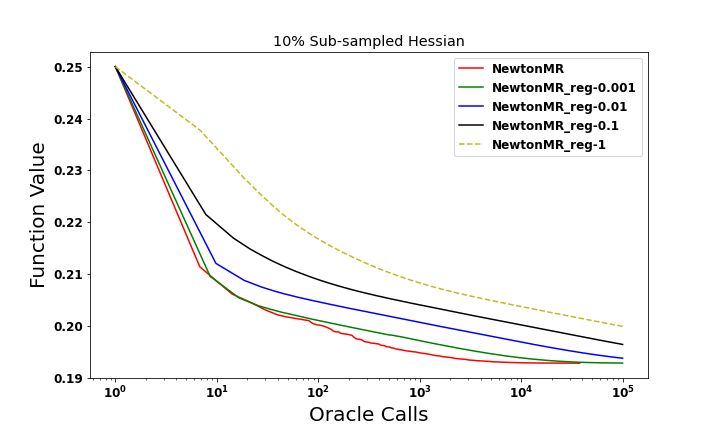} 
\end{subfigure}
\begin{subfigure}{0.33\textwidth}
\includegraphics[width=1\linewidth]{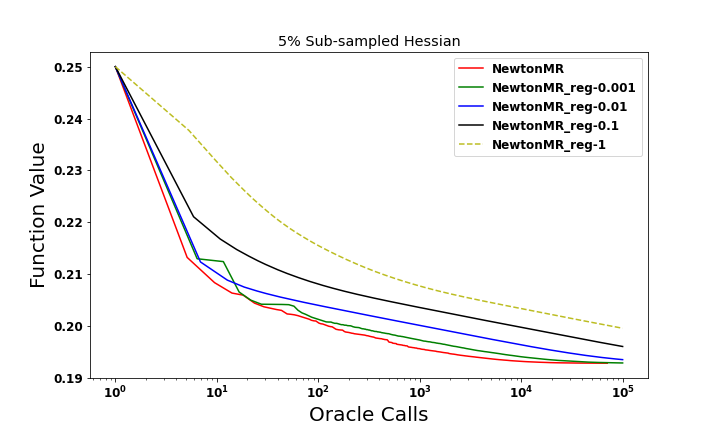} 
\end{subfigure}
\caption{Nonlinear least square problem from \cref{sec:exp:binary} with CIFAR10 dataset \cite{krizhevsky2009learning}. Newton-MR with different magnitude of regularization/damping of (Left) the full Hessian, (Middle) 10\% sub-sampled Hessian and (Right) 5\% sub-sampled Hessian. Here, \texttt{NewtonMR\_reg-$\zeta$} refers to the Newton-MR algorithm \cref{alg:NewtonMR} with regularized/damped Hessian, i.e., $\tilde{\f{H}}_k = \inH_k + \zeta \f{I}$, where $\zeta = 0, 0.001, 0.01, 0.1, 1$. As it can be seen, almost always, such a regularization/damping results in poor empirical performance.\label{fig:damping}}
\end{figure}

As mentioned earlier, our primary focus in this paper is to design a convergent algorithm for nonconvex problems that requires minimal algorithmic modification, relative to similar Newton-type methods in the literature. In nonconvex settings, many Newton-type methods often consider some modifications to the Hessian, e.g. Levenberg-Marquardt  \cite{levenberg1944method,marquardt1963algorithm}, Newton-CG \cite{royer2020newton,xu2020newton}, and Trust-region-CG \cite{curtis2021trust}. In particular, the regularized/damped Hessian in the form of $\tilde{\f{H}}_k = \f{H}_k + \zeta \f{I}$, for some $\zeta > 0$, is used in \cite{royer2020newton,xu2020newton,curtis2021trust} within the Capped-CG procedure. Not only does this particular modification allow for the extraction of sufficiently negative curvature direction, but it also provides a means for obtaining state-of-the-art convergent analysis. However, as we now show, this damping of the Hessian spectrum can adversely affect the practical performance of the algorithm. \cref{fig:damping} depicts the results of the nonlinear least square problem from \cref{sec:exp:binary} with CIFAR10 dataset \cite{krizhevsky2009learning}. We consider different magnitudes of damping for the Hessian. As it can be seen, almost always, such a regularization/damping results in poor empirical performance. Intuitively, this is largely due to the lessening the contribution of the eigenvectors corresponding to the small eigenvalues of the Hessian in forming the update direction.

\end{document}